\documentclass[a4paper,11pt]{amsart}
\usepackage{amssymb}
\usepackage{amscd}

\newtheorem{theorem}{Theorem}[section]
\newtheorem{lemma}[theorem]{Lemma}

\newtheorem{proposition}[theorem]{Proposition}
\newtheorem{corollary}[theorem]{Corollary}

\theoremstyle{definition}

\theoremstyle{remark}

\usepackage{graphicx}
\usepackage{amsmath}
\usepackage{amsfonts}
\usepackage{amssymb}
\newcommand{\be}{\begin{equation}}

\newcommand{\ee}{\end{equation}}

\newcommand{\g}{g^o}

\newcommand{\Om}{\Omega}

\newcommand{\om}{\omega}

\newcommand{\cB}{{\mathcal B}}
\newcommand{\cN}{\mbox{$\mathcal{N}$}}

\newcommand{\fD}{\mbox{$
\begin{picture}(9,8)(1.6,0.15)
\put(1,0.2){\mbox{$ D \hspace{-7.8pt} /$}}
\end{picture}$}}

\newcommand{\boldnabla}{\mbox{\boldmath$ \nabla$}}
\newcommand{\nda}{\boldnabla}

\newcommand{\afD}
{\mbox{$
\begin{picture}(9,8)(1.6,0.15)
\put(1,0.2){\mbox{\boldmath$ D \hspace{-7.8pt} /$}}
\end{picture}$}}

\newcommand{\al}{\mbox{\boldmath$\Delta$}}

\newcommand{\afl}{\mbox{$
\begin{picture}(9,8)(1.6,0.15)
\put(1,0.2){\mbox{$ {\mbox{\boldmath$\Delta$}} \hspace{-7.8pt} /$}}
\end{picture}$}}

\newcommand{\si}{\sigma}

\newcommand{\R}{{\bf R }}

\newcommand{\ba}{\begin{array}}

\newcommand{\ea}{\end{array}}

\newcommand{\beq}{\begin{eqnarray}}

\newcommand{\eeq}{\end{eqnarray}}

\newtheorem{lm}{lemma}

\newtheorem{thee}{theorem}

\newtheorem{proo}{proposition}

\newtheorem{co}{corollary}

\newtheorem{rem}{remark}

\newtheorem{deff}{definition}

\newcommand{\bd}{\begin{deff}}

\newcommand{\ed}{\end{deff}}

\newcommand{\bl}{\begin{lm}}

\newcommand{\el}{\end{lm}}

\newcommand{\bp}{\begin{proo}}

\newcommand{\ep}{\end{proo}}

\newcommand{\bt}{\begin{thee}}

\newcommand{\et}{\end{thee}}

\newcommand{\bc}{\begin{co}}

\newcommand{\ec}{\end{co}}

\newcommand{\brm}{\begin{rem}}

\newcommand{\erm}{\end{rem}}

\usepackage{t1enc}
\def\frak{\mathfrak}
\def\Bbb{\mathbb}
\def\Cal{\mathcal}
\newcommand{\cL}{{\Cal L}}
\newcommand{\bS}{\mathbb{S}}

\newcommand{\newc}{\newcommand}

\let\ccdot\cdot
\def\cdot{\hbox to 2.5pt{\hss$\ccdot$\hss}}

\newc{\aR}{\mbox{\boldmath{$ R$}}}
\newc{\aS}{\mbox{\boldmath{$ S$}}}
\newc{\aT}{\mbox{\boldmath{$ T$}}}
\newc{\aW}{\mbox{\boldmath{$ W$}}}
\newcommand{\aX}{\mbox{\boldmath{$ X$}}}

\newc{\aK}{\mbox{\boldmath{$ K$}}}
\newc{\aL}{\mbox{\boldmath{$ L$}}}

\newcommand{\ce}{{\Cal E}}

\newcommand{\cq}{{\Cal Q}}

\newcommand{\ct}{{\Cal T}}

\newcommand{\bV}{{\Bbb V}}

\newcommand{\bD}{\mathbb{D}}

\usepackage{amssymb}
\usepackage{amscd}

\newcommand{\nd}{\nabla}

\newcommand{\Rho}{P}
\newcommand{\Up}{\Upsilon}

\newcommand{\End}{\operatorname{End}}

\newcommand{\Ric}{\operatorname{Ric}}
\newcommand{\Sc}{\operatorname{Sc}}
\newcommand{\Proj}{\operatorname{Proj}}

\newcommand{\ul}[1]{\underline{#1}}

\newcommand{\bs}[1]{\boldsymbol{#1}}

\def\endrk{\hbox{$|\!\!|\!\!|\!\!|\!\!|\!\!|\!\!|$}}
\newcommand{\cT}{{\mathcal T}}

\let\hash=\sharp  

\newcommand{\fb}  {\mbox{$                                      
\begin{picture}(11,8)(0.5,0.15)
\put(1,0.2){\mbox{$ \Box \hspace{-7.8pt} /$}}
\end{picture}$}} 

\def\endrk{\hbox{$|\!\!|\!\!|\!\!|\!\!|\!\!|\!\!|$}}

\let\i=\iota

\newcommand{\hook}{\raisebox{-0.35ex}{\makebox[0.6em][r]
{\scriptsize $-$}}\hspace{-0.15em}\raisebox{0.25ex}{\makebox[0.4em][l]{\tiny
 $|$}}}

\newcommand{\aM}{\tilde{M}}

\newcommand{\act}{\mbox{\boldmath{$\ct$}}}

\newcommand{\nn}[1]{(\ref{#1})}

\newcommand{\X}{\mbox{\boldmath{$ X$}}}
\newcommand{\sX}{\mbox{\scriptsize\boldmath{$X$}}}        
\newcommand{\h}{\mbox{\boldmath{$ h$}}}
\newcommand{\bg}{\mbox{\boldmath{$ g$}}}
\newcommand{\sbg}{\mbox{\boldmath{\scriptsize$ g$}}}

\renewcommand{\S}{\Sigma}

\newcommand{\J}{{\mbox{\sf J}}}
\newc{\obstrn}[2]{B^{#1}_{#2}}

\newcommand{\rpl}                         
{\mbox{$
\begin{picture}(12.7,8)(-.5,-1)
\put(0,0.2){$+$}
\put(4.2,2.8){\oval(8,8)[r]}
\end{picture}$}}

\newcommand{\lpl}                         
{\mbox{$
\begin{picture}(12.7,8)(-.5,-1)
\put(2,0.2){$+$}
\put(6.2,2.8){\oval(8,8)[l]}
\end{picture}$}}

\usepackage{ifthen}

\newc{\tensor}[1]{#1}
\newc{\Mvariable}[1]{\mbox{#1}}
\newc{\down}[1]{{}_{#1}}
\newc{\up}[1]{{}^{#1}}

\newc{\JulyStrut}{\rule{0mm}{6mm}}
\newc{\midtenPan}{\mbox{\sf S}}
\newc{\midten}{\mbox{\sf T}}
\newc{\midtenEi}{\mbox{\sf U}}
\newc{\ATen}{\mbox{\sf E}}
\newc{\BTen}{\mbox{\sf F}}
\newc{\CTen}{\mbox{\sf G}}

\renewcommand{\P}{\mbox{$\Bbb P$}}

\def\sideremark#1{\ifvmode\leavevmode\fi\vadjust{\vbox to0pt{\vss
 \hbox to 0pt{\hskip\hsize\hskip1em
 \vbox{\hsize3cm\tiny\raggedright\pretolerance10000
 \noindent #1\hfill}\hss}\vbox to8pt{\vfil}\vss}}}%
                        
\numberwithin{equation}{section}

\begin{document}
\renewcommand{\today}{}
\title{ Almost Einstein and
Poincar\'e-Einstein manifolds in Riemannian signature}
\author{A. Rod Gover}

\address{Department of Mathematics\\
  The University of Auckland\\
  Private Bag 92019\\
  Auckland 1\\
  New Zealand} \email{gover@math.auckland.ac.nz}

\vspace{10pt}

\renewcommand{\arraystretch}{1}

\begin{abstract} 
 An almost Einstein manifold satisfies equations which are a slight
  weakening of the Einstein equations; Einstein metrics,
  Poincar\'e-Einstein metrics, and compactifications of certain
  Ricci-flat asymptotically locally Euclidean structures are special
  cases. The governing equation is a conformally invariant
  overdetermined PDE on a function. Away from the zeros of this the
  almost Einstein structure is Einstein, while the zero set gives a
  scale singularity set which may be viewed as a conformal infinity for
  the Einstein metric.  In this article we give a classification of
  the possible scale singularity spaces and derive geometric results
  which explicitly relate the intrinsic conformal geometry of these to
  the conformal structure of the ambient almost Einstein manifold.
  Classes of examples are constructed.  A compatible generalisation of
  the constant scalar curvature condition is also developed. This
  includes almost Einstein as a special case, and when its curvature
  is suitably negative, is closely linked to the notion of an
  asymptotically hyperbolic structure.
\end{abstract}

\maketitle
\renewcommand{\arraystretch}{1.5}

\pagestyle{myheadings} \markboth{Gover}{Almost Einstein and
Poincar\'e-Einstein manifolds in Riemannian signature}

\thanks{The author gratefully
  acknowledges support from the Royal Society of New Zealand via
  Marsden Grant no.\ 06-UOA-029}

\maketitle
\section{Introduction}

A metric is said to be Einstein if its Ricci curvature is proportional
to the metric \cite{Besse}. Despite a long history of intense interest
in the Einstein equations many mysteries remain.  In high dimensions
it is not known if there are any obstructions to the existence of
Einstein metric. There are 3-manifolds and 4-manifolds which do not
admit Einstein metrics and the situation is especially delicate in the
latter case, see \cite{leBrun0404359} for an overview of some recent
progress.  Here we consider a specific weakening of the Einstein
condition. By its nature this provides an alternative route to
studying Einstein metrics but, beyond this, there are several points
which indicate that it may be a useful structure in its own right. On
the one hand the weakening is very slight, in a sense that will soon
be clear. On the other it allows in some interesting cases: at least
some manifolds satisfying these equations do not admit Einstein
metrics, which suggests a role as a uniformisation type condition; it
includes in a natural way Poincar\'e-Einstein structures and
conformally compact Ricci-flat asymptotically locally Euclidean (ALE)
spaces, and so Einstein metrics, Poincar\'e-Einstein structures and these
ALE manifolds are special cases of a uniform generalising structure.

 On a Riemannian manifold $(M^d,g)$
($d\geq 3$ here and throughout) the Schouten tensor $P$ (or $P^{g}$)
is a trace adjustment of the Ricci tensor given by
$$
\Ric^g=(d-2)P^g+J^g g
$$ where $J^g$ is the metric trace of $P^g$. Thus a metric is Einstein if and
only if the trace-free part of $P^g$ is zero.  We will say
that $(M,g,s)$ is a {\em directed almost Einstein} structure if $s\in
C^\infty (M)$ is a non-trivial solution to the equation
\begin{equation}\label{prim}
A(g,s)=0 \quad \mbox{where} \quad A(g,s) :=trace-free(\nabla^g\nabla^g
s + sP^g).
\end{equation} 
Here $\nabla^g$ is the Levi-Civita connection for $g$, and the
``trace-free'' means the trace free part with respect to taking a
metric trace. This is a generalisation of the Einstein condition; we
will see shortly that, on the open set where $s$ is non-vanishing,
$\g:=s^{-2}g$ is Einstein; here Einstein is forced as a consistency
condition for a solution to \nn{prim}.  On the other hand if $g$ is
Einstein then \nn{prim} holds with $s=1$. Any attempt to understand
the nature and extent of this generalisation should include a
description of the possible local structures of the {\em scale
singularity set}, that is the set $\Sigma$ where $s$ is zero (and
where $\g=s^{-2}g$ is undefined). The main results in this article are
some answers to this question and the development of a conformal
theory to relate, quite directly, the intrinsic geometric structure of
the singularity space $\S$ to the ambient structure.  If $s$ solves
\nn{prim} then so does $-s$, and where $s$ is non-vanishing these
solutions determine the same Einstein metric. We shall say that a
manifold $(M,g)$ is {\em almost Einstein} if it admits a covering such
that on each open set $U$ of the cover we have that $(U,g,s_U)$ is
directed almost Einstein and on overlaps $U\cap V$ we have either
$s_U=s_V$ or $s_U=-s_V$. Although there exist almost Einstein spaces
which are not directed \cite{GoLeitprog}, to simplify the exposition
we shall assume here that almost Einstein (AE) manifolds are
directed. (So usually we omit the term ``directed'' but sometimes it
is included in Theorems for emphasis.) In any case the results apply
locally on almost Einstein manifolds which are not directed.

On an Einstein manifold $(M,g)$ the Bianchi identity implies that the
scalar curvature $\Sc^g$ (i.e.\ the metric trace of $\Ric$) is
constant. Thus simply requiring a metric to be scalar constant is another 
weakening of the Einstein condition. On compact, connected oriented
smooth Riemannian manifolds this may be achieved conformally: this is
the outcome of the solution to the ``Yamabe problem'' due to Yamabe,
Trudinger, Aubin and Schoen \cite{Y,T,A,S}.  Just as almost Einstein
generalises the Einstein condition, there is an corresponding
weakening of the constant scalar curvature condition as follows.  We
will say that $(M,g,s)$ is a {\em directed almost scalar constant}
structure if $s\in C^\infty (M)$ is a non-trivial solution to the
equation $S(g,s)=constant$ where
\begin{equation}\label{ascE}
S(g,s)= \frac{2}{d}s(J^g-\Delta^g)s - |ds|^2_g~.  
\end{equation} 
Away from the zero set (which again we denote by $\Sigma$) of $s$ we
have $S(g,s)=\Sc^{\g}/d(d-1)$ where $\g:=s^{-2}g$. In particular, off
$\Sigma$, $S(g,s)$ is constant if and only if $\Sc^{\g}$ is constant.
The normalisation is so that if $\g$ is the metric of a space form
then $S(g,s)$ is exactly the sectional curvature.  We shall say that a
manifold $(M,g)$ is {\em almost Scalar constant} (ASC) if it is equipped
with a covering such that on each open set $U$ of the cover we have
that $(U,g,s_U)$ is directed almost scalar constant, and on overlaps $U\cap
V$ we have either $s_U=s_V$ or $s_U=-s_V$. In fact, in line with our 
assumptions above and unless otherwise mentioned explicitly, 
 we shall assume below that any ASC structure is directed.

As suggested above, closely related to these notions are certain
classes of so called conformally compact manfolds that have recently
been of considerable interest.  We recall how these manifolds are
usually described.  Let $M^{d}$ be a compact smooth manifold with
boundary $\Sigma=\partial M$. A metric $g^o$ on the interior $M^+$ of
$M$ is said to be conformally compact if it extends (with some
specified regularity) to $\Sigma$ by $g=s^2g^o$ where $g$ is
non-degenerate up to the boundary, and $s$ is a non-negative defining
function for the boundary (i.e.\ $\Sigma$ is the zero set for $s$, and
$d s$ is non-vanishing along $M$).  In this situation the metric $g^o$
is complete and the restriction of $g$ to $T\Sigma$ in $TM|_\Sigma$
determines a conformal structure that is independent of the choice of
defining function $s$; then $\S$ with this conformal structure is
termed the conformal infinity of $M^+$. (This notion had its origins
in the work of Newman and Penrose, see the introduction of
\cite{LeBrunH} for a brief review.)  If the defining function is
chosen so that $|ds|^2_g=1$ along $M$ then the sectional curvatures
tend to $-1$ at infinity and the structure is said to be
asymptotically hyperbolic (AH) (see \cite{M-hodge} where there is a
detailed treatment of the Hodge cohomology of these structures and
related spectral theory).  The model is the Poincar\'e hyperbolic ball
and thus the corresponding metrics are sometimes called Poincar\'e
metrics.  Generalising the hyperbolic ball in another way, one may
suppose that the interior conformally compact metric $g^o$ is Einstein
with the normalisation $\Ric (g^o)=-n g^o$, where $n=d-1$, and in this
case the structure is said to be Poincar\'e-Einstein (PE); in fact PE
manifolds are necessarily asymptotically hyperbolic.  Such structures
have been studied intensively recently in relation to the proposed
AdS/CFT correspondence of Maldacena \cite{Mal,Witten}, related
fundamental geometric questions
\cite{Albin,And,Biquard,CQYrenromvol,GL,GrZ,lee,MazPac} and through
connections to the ambient metric of Fefferman-Graham
\cite{FGast,FGnew}.

For simplicity of exposition we shall restrict our attention to smooth
AE and ASC structures $(M^{d},g,s)$; that is $(M,g)$ is a smooth
Riemannian manifold of dimension $d\geq 3$ and $s\in C^\infty(M)$
satisfies either \nn{prim} (the AE case) or \nn{ascE} (for ASC).  Let
us write $M^\pm$ for the open subset of $M$ on which $s$ is positive
or, respectively, negative and, as above, $\S$ for the scale
singularity set. The first main results (proved in Section \ref{al}) are the following classifications
for the possible submanifold structures of $\Sigma$.
\begin{theorem} \label{mainsc} Let $(M,^dg,s)$ be a directed almost scalar
  constant structure with $g$ positive definite and $M$ connected.  If
  $S(g,s) > 0$ then $s$ is nowhere vanishing and $(M,\g)$ has constant
  scalar curvature $d(d-1)S(g,s)$. If $S(g,s) < 0$ then $s$ is
  non-vanishing on an open dense set and $\S$ is either empty or else
  is a hypersurface; On $M\setminus \S$, $\Sc^{\g}$  is constant and equals 
$d(d-1)S(g,s)$.  Suppose $M$ is closed (i.e.\ compact without boundary)
  with $S(g,s)<0$ and $\S\neq \emptyset$. A constant rescaling of $s$
  normalises $S(g,s)$ to $-1$, and then $(M\setminus M^-)$ is a finite
  union of connected AH manifolds. Similar for $(M\setminus M^+)$.
\end{theorem} 
\noindent By hypersurface we mean a submanifold of codimension 1 which
may include boundary components. In the following we will say that an
ASC structure is scalar positive, scalar flat, or scalar negative if,
respectively, $S(g,s)$ is positive, zero, or negative.

It seems that almost Einstein manifolds, in the generality we describe
here, were introduced in \cite{GoPrague} and it was observed there
that PE manifolds are a special case; this was explained in detail in
\cite{GoMGC}. Here, among other things, we see that PE manifolds arise
automatically in the scalar negative (i.e.\ $S(g,s)<0$) case.  
\begin{theorem} \label{maina} Let $(M,g,s)$ be a directed  almost Einstein 
structure with $g$ positive definite and $M$ connected. Then $s$ is
non-vanishing on an open dense set and $(M,g,s)$ is almost scalar
constant.  Writing $\S$ for the scale singularity set, on $M\setminus
\S$, $\g$ is Einstein with scalar curvature $d(d-1)S(g,s)$. There are
three cases: \\ $\bullet$ If $S(g,s)>0$ then the scale singularity set
$\S$ is empty. \\ $\bullet$ If $S(g,s)=0$ then $\S$ is either empty or
otherwise consists of isolated points and these points are critical
points of the function $s$; in this case for each $p\in M$ with
$s(p)=0$, the metric $\g$ is asymptotically locally Euclidean (ALE)
near $p$ and the Weyl, Cotton, and Bach curvatures vanish at $p$. \\
$\bullet$ If $S(g,s)<0$ then $\S$ is either empty or else is a totally
umbillic hypersurface.  In particular on a closed $S(g,s)=-1$ almost
Einstein manifold $(M\setminus M^-)$ is a finite union of connected
Poincar\'e-Einstein manifolds. Similar for $(M\setminus M^+)$.
\end{theorem}
The Cotton and Bach curvatures are defined in, respectively, \nn{cot}
and \nn{Bdef} below.
Using compactness, the last statement is an easy consequence of
 Proposition \ref{peae}. That AE implies ASC is part of Theorem
 \ref{key}.  Given this several parts of the Theorem are immediate
 from Theorem \ref{mainsc} above. The remaining parts of the Theorem
 summarise Theorem \ref{classthm}, Proposition \ref{nullp},
 Proposition \ref{aEu}, and parts of Proposition \ref{aEsum} and
 Corollary \ref{weakgen}. We shall say that the ALE structures arising as
 here are conformally {\em conformally compact} because of the obvious
 link the term as used above.

The equation \nn{ascE} is conformally covariant in the sense that for
any $\om\in C^\infty(M)$ we have $S(g,s)=
S(e^{2\om}g,e^{\om}s)$. Similarly for \nn{prim} we have $e^\om A(g,s)=
A(e^{2\om}g,e^{\om}s)$ and so if $(M,g,s)$ is almost Einstein then so
is $(M,e^{2\om}g,e^{w}s)$. Evidently the notions of ASC and AE
structure pass to the conformal geometry by taking a quotient of the
space of all such structures by the equivalence relation $(M,g,s)\sim
(M,e^{2\om}g,e^{w}s)$. This is the point of view we wish to take,
throughout $g$ is to be viewed as simply a representative of its
conformal class.  (This shows that we should really view the function
$s$ as corresponding, via the density bundle trivialisation afforded
by the metric $g$, to a conformal density $\si$ of weight 1 on the
conformal manifold $(M,[g])$, and $A$ as a 2-tensor taking values in
this density bundle. We shall postpone this move until Section
\ref{al}.)  The conformal equivalence class of $(g,s)$ (under $
(g,s)\sim (e^{2\om}g,e^{\om} s)$) is a structure which generalises the
notion of a metric.  This suggests a definition which is convenient
for our discussions. A Riemannian manifold equipped with the conformal
equivalence class (in this sense) of $(g,s)$, and where $s$ is nowhere
vanishing on an open dense set, is a well defined structure that we shall term
an {\em almost Riemannian} manifold. Of course the zero set of $s$ is
conformally invariant and so is a preferred set $\Sigma\subset M$.
 An almost Riemannian structure with
$\S=\emptyset$ is simply  a Riemannian manifold. 
Note that in the other cases
$S(g,s)$ smoothly extends to $M$ the natural scalar $Sc^{\g}/{n(n+1)}$
which is only defined on $M\setminus \Sigma$.  Similarly $A(g,s)$
smoothly extends $s P^{\g}_0$, where $P^{\g}_0$ is the trace free part
of $P^{\g}$. Thus even though the metric $\g=s^{-2}g$ is not defined
along $\S$, nevertheless $A(g,s)$ and $S(g,s)$ are defined globally
(at least if we view $A(g,s)$ as representing a density valued tensor)
and it is natural to think of these as curvature quantities on almost
Riemannian structures. It turns out that AE
manifolds, and also the cases of ASC manifolds covered in Theorem
\ref{mainsc}, are necessarily almost Riemannian.

The structures we consider here have an elegant and calculationally
effective formulation in terms of conformal tractor calculus.  On
Riemannian manifolds the metric canonically determines a connection on
the tangent bundle, the Levi-Civita connection. On conformal
structures we lose this but there is a canonical conformally invariant
connection $\nabla^{\cT}$ on the (standard conformal) {\em tractor
  bundle} $\cT$, as described in the next section. On $(M^d,[g])$ this
is a rank $(d+2)$ bundle that contains a conformal twisting of the
tangent bundle as a subquotient.  The bundle $\cT$ also has a
(conformally invariant) tractor metric $h$, of signature $(d+1,1)$,
that is preserved by $\nabla^\cT$. On a given conformal structure we
may ask if there is parallel section of $\cT$; that is a secction $I$
of $\cT$ satisfying $\nabla^\ct I=0$.  In fact, as we see below (following \cite{BEG}), this
equation is simply a prolongation of \nn{prim}. In particular, on any
open set, solving $\nabla^\ct I=0$ is equivalent to solving \nn{prim}
and there is an explicit 1-1 relationship between solutions. (We shall
write $s_I$ for the solution of \nn{prim} given by a parallel tractor
$I$.)  Thus an almost Einstein structure is a triple $(M,[g],I)$ where
$I$ is parallel for the standard tractor connection determined by the
conformal structure $[g]$. Since the tractor connection preserves the
metric $h$, the length (squared) of $I$, which we denote by the
shorthand $|I|^2:=h(I,I)$, is constant on connected AE manifolds (and
we henceforth assume $M$ is connected).  In fact $S(g,s_I)= -|I|^2 $. There is 
a hgeneralising result for ASC manifolds, see Proposition \ref{ascP}.

The geometric study of PE manifolds has been driven by a desire to
relate the conformal geometry of the conformal infinity to the metric
geometry on the interior. We may obviously extend this programme to
the scalar negative (i.e.\ $S(g,s)<0$) almost Einstein structures.  As
indicated above, this is a core aim here and in our treatment
(Sections \ref{MvS} and \ref{FGsec}) the tractor structures play a key
role.  The first key result is Theorem \ref{umbiltr} which shows, for
example, that $\S$ satisfies a conformal analogue of the Riemannian
totally geodesic condition: the intrinsic tractor connection of
$(\S,[g_\S])$ exactly agrees with a restriction of the ambient tractor
connection. In fact the results are stronger. Summarising part of
Theorem \ref{umbiltr} with Corollary \ref{WW}, along the scale
singularity set $\Sigma$ of a scalar negative AE structure we also
have the following:
\begin{theorem} \label{curvTr} 
$$
\Omega(u,v)=\Omega^\Sigma(u,v) \quad \mbox{ along } \Sigma
$$
where $u,v\in \Gamma(T\Sigma)$. In dimensions $d\neq 4$ we have the stronger result
$$
\Omega(\cdot,\cdot)=\Omega^\Sigma(\cdot,\cdot)  \quad \mbox{ along } \Sigma ,
$$ where here, by trivial extension, we view $\Omega^\Sigma$ as a section of 
$\Lambda^2 T^*M\otimes \End \cT$.
While in dimensions $d\geq 6$ we also have 
$$
(d-5)W|_\S=(d-4)W^\S,
$$ where $W$ is the prolonged conformal curvature quantity \nn{Wform}
and again a trivial extension is involved.
\end{theorem}
Here $\Omega$ is the curvature of the tractor connection for $(M,[g])$
while $\Omega^\Sigma$ is the curvature of the tractor connection for
the intrinsic conformal structure of $\Sigma$. $W$ is the natural
conformally invariant tractor field equivalent (in dimensions $d\neq
4$) to the curvature of the Fefferman-Graham (ambient) metric over
$(M,[g])$, while $W^\S$ is the same for $(\S,[g_\S])$.  In Section
\ref{FGsec} Theorem \ref{FGflat} we also show that the
Fefferman-Graham (obstruction) tensor must vanish on the scale
singularity hypersurface of an almost Einstein structure. An
alternative direct proof that $\S$ is Bach-flat, when $n=4$, is given
in Corollary \ref{Bflat}.  A key tool derived in Section \ref{FGsec}
is Theorem \ref{amb} which constructs a Fefferman-Graham ambient
metric, formally to all orders, for the even dimensional conformal
structure of a scale singularity set; this construction was heavily
influenced by the model in Section \ref{model}. An important and
central aspect of the works \cite{FGast} and \cite{FGnew} is the
direct relationship between the Fefferman-Graham (ambient) metric for
conformal manifolds $(\S,[g_\S])$ and suitably {\em even} smooth
formal Poincar\'e-Einstein metrics, with $(\S,[g_\S])$ as the
conformal infinity (see especially \cite[Section 4]{FGnew}); in
Section \ref{dou} there is some discussion of the meaning of even in
this context. Here, in contrast, we work in one higher dimension and
exploit the use of the ambient metric for the Poincar\'e-Einstein (or
AE) space $M$ itself as tool for studying the boundary $\S$ (or scale
singularity set); in this case we may work globally on $M$ and with
not necessarily even PE (or AE) metrics.

In Section \ref{extS} we describe equations controlling (at least
partially) the conformal curvature of almost Einstein
structures. Importantly these are given in a way that should be
suitable for setting up boundary problems along $\S$ based around the
conformal curvature quantities. For example in Proposition \ref{d4} we
observe that in this sense the Yang-Mills equations, applied to the
tractor curvature, give the natural conformal equations for
4-dimensional almost Einstein structures. The anologue for higher even
dimensions is given in Proposition \ref{conv}. In all dimensions we
have the following result.
\begin{theorem}\label{Weqn}
Let $(M^d,[g],I)$ be an almost Einstein manifold then
$$
I^A\fD_A W=0.
$$ 
\end{theorem}
The operator $I^A \fD_A$ has the form $\si \Delta +lower~order~terms$
except along $\S$ in dimensions $d\neq 6$. The statement here is
mainly interesting in dimensions $d\geq 5$ and is a part of Theorem
\ref{ext}.  Since for $d\geq 6$ we have $(d-5)W|_\S=(d-4)W^\S$, for
Poincar\'e-Einstein manifolds (and more generally scalar negative AE
structures) the Theorem suggests a Dirichlet type problem with the
conformal curvature $W^\S$ of $\S$ as the boundary (hypersurface)
data.  The operator $I^A \fD_A$ is well defined on almost Einstein
manifolds and is linked to the scattering picture of \cite{GrZ} as
outlined in Corollary \ref{interp}.

As mentioned, almost Einstein structures provide a generalisation of
the notions of Einstein, Poincar\'e-Einstein and certain conformally
compact ALE metrics. Aside from providing a new and uniform
perspective on these specialisations, the AE structures
provide a natural uniformisation type problem. We may ask for
example whether any closed smooth manifold admits an almost Einstein
structure. While it is by now a classical result \cite{Besse} that the
sphere products $\bS^1\times \bS^2$ and $\bS^1\times \bS^3$ do not
admit Einstein metrics it is shown in \cite{GoLeitprog} that these
both admit almost Einstein structures; in fact we construct these
explicitly as part of a general construction of closed manifolds with
almost Einstein structures. In this article we make just a small
discussion of examples in Section \ref{exsS}. This includes the
conformal sphere as the key model. It admits all scalar types of
almost Einstein structure and has a central role in the construction
of other examples in \cite{GoLeitprog}.  (In fact the standard
conformal structure on the sphere admits a continuous curve of AE
structures which includes the standard sphere metric, the
Euclidean metric pulled back by stereographic projection as well as negative
$S(g,s)$ AE structures with $\g$  hyperbolic off the singularity set. See
Corollary \ref{betw} and the final comments in Section \ref{model}.)
We conclude in Section \ref{dou} with a discussion of examples found
by a doubling construction. Non-Einstein almost Einstein metrics turn
up in the constructions and classifications by Derdzinski and Maschler
of K\"ahler metrics which are ``almost everywhere'' conformal to
Einstein by a non-constant recaling factor, see e.g. \cite{DML,DM} and
references therein. Some of their examples were inspired by
constructions known for some time, such as \cite{D,P}.  Examples of
non-Einstein $S([g],I)=0$ AE structures are disussed in \cite{KR}.

 It should also be pointed out that many of the techniques and results
we develop apply in other signatures. However there are also
fundamental differences in the case of non-Riemannian signature and so
we confine the study to the positive definite setting.

Conversations with Michael Eastwood, Robin Graham, Felipe Leitner, and
Paul-Andi Nagy have been much appreciated.  It should pointed out that
 the existence of AE structures which are not directed was
observed in the joint work \cite{GoLeitprog} with Leitner and this
influenced the presentation here.

\section{Almost Einstein structures and conformal tractor calculus}\label{al}

As above let $M$ be a smooth manifold, of dimension $d\geq 3$,
equipped with a Riemannian metric $g_{ab}$. Here and throughout we
employ Penrose's abstract index notation. We write $\ce^a$ to denote
the space of smooth sections of the tangent bundle $TM$ on $M$, and $\ce_a$
for the space of smooth sections of the cotangent bundle $T^*M$.  (In fact we
will often use the same symbols for the 
bundles themselves. Occasionally, to avoid any confusion, we write
$\Gamma(\mathcal{B})$ to mean the space of sections of a bundle
$\mathcal{B}$.) We write $\ce$ for the space of smooth functions and
all tensors considered will be assumed smooth without further comment.
An index which appears twice, once raised and once lowered, indicates
a contraction.  The metric $g_{ab}$ and its inverse $g^{ab}$ enable
the identification of $\ce^a$ and $\ce_a$ and we indicate this by
raising and lowering indices in the usual way.

With $\nabla_a$ denoting the Levi-Civita connection for $g_{ab}$, and
using that this is torsion free, the Riemann curvature tensor
$R_{ab}{}^{c}{}_d$ is given by
$$(\nabla_a\nabla_b-\nabla_b\nabla_a)V^c=R_{ab}{}^{c}{}_d V^d 
\quad\text{
where} \quad \ V^c\in \ce^c.$$
This can be decomposed into the totally trace-free {\em Weyl curvature}
$C_{abcd}$ and the symmetric {\em
Schouten tensor} $\Rho_{ab}$ according to
\begin{equation}\label{Rsplit}
R_{abcd}=C_{abcd}+2g_{c[a}\Rho_{b]d}+2g_{d[b}\Rho_{a]c},
\end{equation}
where $[\cdots]$ indicates  antisymmetrisation over the enclosed
indices. 
Thus $\Rho_{ab}$ is a trace modification of the Ricci tensor 
${\rm Ric}_{ab}=R_{ca}{}^c{}_b$:
$$
\Ric_{ab}=(n-2)\Rho_{ab}+ J g_{ab}, \quad \quad J:=\Rho^a_{~a}.
$$ 
In denoting such curvature quantities we may write e.g.\ $\Ric^g$ or
simply $\Ric$ depending on whether there is a need to emphasise the
metric involved.  Also abstract indices will be displayed or suppressed
as required for clarity.

Under a conformal rescaling of the metric 
$$
g\mapsto \g=s^{-2}g, 
$$ with $s\in \ce$ non-vanishing, the Weyl tensor $C_{ab}{}^c{}_d$ is
is unchanged (and so we say the Weyl tensor is conformally invariant)
whereas the Schouten tensor transforms according to
\begin{equation}\label{Pt}
 P^{\g}_{ab}= P^g_{ab} +s^{-1}\nd_a\nd_b s -\frac{1}{2}g^{cd}s^{-2}(\nd_c
s)(\nd_d s)g_{ab}.
\end{equation}
Taking, via $\g$, a trace of this we obtain 
$$
J^{\g}=s^2 J^g -s\Delta s - \frac{d}{2}|ds|_g^2,
$$
where the $\Delta$ is the ``positive energy'' Laplacian.
Note that the right hand side is $\frac{d}{2}S(g,s)$, with $S(g,s)$ as
defined in \nn{ascE}. Clearly this is well defined for smooth $s$ even
if $s$ may be  zero at some points. On the other hand the right hand side above
(and hence $S(g,s)$) is clearly invariant under the conformal
transformation $(g,s)\mapsto (e^{2\om}g,e^{\om}s)$: this is true away
from the zeros of $s$ since there $J^{\g}$ depends only on the 2-jet
of $\g=s^{-2}g$, but the explicit conformal transformation of the
right hand side is evidently polynomial in $e^\om$ and its 2-jet.

Let us digress to prove Theorem \ref{mainsc} since it illustrates how
an almost Riemannian structure may arise immediately from a formula
polynomial in the jets of $s$.\\
\noindent{\bf Proof of Theorem \ref{mainsc}:} Under a dilation
$g\mapsto \mu g$ ($\mu\in \mathbb{R}_+$) we have $S(g,s)\mapsto
\mu^{-1} S(g,s)$, so to prove the Theorem we may consider just the
cases $S(g,s)=1$ and $S(g,s)=-1$.  Suppose that $S(g,s)=1$ then if
$p\in M$ were to be a point where $s_p=0$ then at $p$ we would have
$1=-|ds|^2_g$ which would be a contradiction.  Suppose that
$S(g,s)=-1$. Then at any point $p\in M$ where $s_p=0$ we have
$|ds|^2_g=1$.  For the last statement of the Theorem assume that $M$
is closed and the scale singularity set $\S$ is not empty.  Then $\S$
is a hypersurface which separates $M$ according to the sign of
$s$. The restriction of $\g$ to the interior of $M\setminus M^-$
(i.e.\ to $M_+$) is conformally compact since the restriction of $g$
to $M\setminus M^-$ extends $s^2\g$ smoothly to the boundary.  Finally
$(M\setminus M^-,g,s)$ is AH since $|ds|^2_g=1$ along $\Sigma$. By
compactness this consists of a finite union of connected AH
components. The same analysis applies to $M\setminus M^+$. \quad
$\Box$. \\
 Note that although constant $S(g,s)$ is a weakening of the
constant scalar curvature condition, the equation \nn{ascE} is quite
restrictive. For example, it is evident that on closed manifolds with
negative Yamabe constant there are no non-trivial solutions with
$S(g,s)$ a non-negative constant.

The tensor $A(g,s)$ defined in the Introduction should be compared to
the trace free part of the right hand side of \nn{Pt} above. Arguing
as for $S(g,s)$ above, or by direct calculation, one finds that under
$(g,s)\mapsto (e^{2\om}g,e^{\om}s)$ we have $A(g,s)\mapsto e^\omega
A(g,s)$ as mentioned earlier. So both the AE condition and the more
general ASC condition are best treated as structures on a conformal
manifold. To obtain a clean treatment it is most efficient to draw in
some standard objects from conformal geometry; for these further
details and background may be found in \cite{CapGoamb,GoPetCMP}.
Clearly we may view a conformal structure on $M$ is a
smooth ray subbundle $\cq\subset S^2T^*M$ whose fibre over $x$
consists of conformally related  metrics at the point
$x$. The principal bundle $\pi:\cq\to M$ has
structure group $\Bbb R_+$, and so each representation ${\Bbb R}_+ \ni
x \mapsto x^{-w/2}\in {\rm End}(\Bbb R)$ induces a natural line bundle
on $ (M,[g])$ that we term the conformal density bundle $\ce[w]$. We
shall write $ \ce[w]$ for the space of sections of this bundle.
 Note
$\ce[w]$ is trivialised by a choice of metric $g$ from the conformal
class, and we write $\nd$ for the connection corresponding to this
trivialisation.  It follows immediately that (the coupled) $ \nd_a$
preserves the conformal metric. (Note on a fixed conformal structure
the conformal densities bundle $\ce[-n]$  may be
identified in an obvious way with appropriate powers of the 1-density
bundle associated to the frame bundle through the representation
$|\det(~)|^{-1}$. See e.g.\ \cite{CapGoamb}. Via this the connection
we defined on $\ce[w]$ agrees with the usual Levi-Civita connection.)

We write $\bg$ for the {\em conformal metric}, that is the
tautological section of $S^2T^*M\otimes \ce[2]$ determined by the
conformal structure. This will be henceforth used to identify $TM$
with $T^*M[2]$ even when we have fixed a metric from the conformal
class. (For example, with these conventions the Laplacian $ \Delta$ is
given by $\Delta=-\bg^{ab}\nd_a\nd_b= -\nd^b\nd_b\,$.)  Although this
is conceptually valuable and significantly simplifies many
calculations, it is, however,  a point where there is potential for confusion.
For example in the below, when we write $J$ or $J^g$ we mean
$\bg^{ab}P_{ab}$ where $P$ is the Schouten tensor for some 
metric $g$. Thus $J$ is a section of $\ce[-2]$ (which depends on
$g$).  

In this picture to study the ASC condition we replace $s\in
\ce$ with a section $\si\in \ce[1]$ in \nn{ascE}  
to obtain
\begin{equation}\label{ascC}
S([g],\si)= \frac{2}{d} \si (J^g-\Delta^g)\si - |\nd \si |^2_{\bg}~,
\end{equation}
where we have written $|\nd\si |^2_{\bg} $ as a brief notation for
$\bg^{-1}(\nabla \si, \nd \si)$.  When the conformal structure is
fixed we shall often denote the quantity displayed by simply
$S(\si)$.  Similarly the conformally invariant version of $A$ is the
2-tensor of conformal weight 1 given by
$$ 
A([g],\si):= trace-free(\nd_a\nd_b \si+P_{ab} \si),
$$
again we may write simply $A(\si)$. 

The $A([g],\si)=0$ equation (i.e.\ \nn{prim}) becomes
\begin{equation}\label{primc}
\nd_a\nd_b \si+P_{ab} \si+\rho \bg_{ab}=0
\end{equation} 
where $\rho$ is an unknown density (in $\ce[-1]$) to accommodate the
trace-part. Here $\nd$ and $P$ are given with respect to some metric
$g$ in the conformal class, but the equation is invariant
under conformal rescaling. 

We may replace \nn{primc} with the
equivalent first order system
$$
\nd_a \si -\mu_a=0, \quad \mbox{and} \quad \nd_a\mu_b+P_{ab}+\bg_{ab}\rho=0,
$$
where $\mu_a\in \ce_a[1]:=\ce_a\otimes \ce[1]$.
Differentiating the second of these and considering two possible  
contractions  yields 
$$
\nd_a \rho -P_{ab}\mu^b=0,
$$ whence we see that the system has closed up linearly. The equation
\nn{primc} is equivalent to a connection and a parallel section for
this; on any open set in $M$, a solution of \nn{primc} is equivalent
to  $I:=(\si,~\mu_a,~\rho)\in \ce[1]\oplus\ce_a[1]\oplus \ce[-1]$ satisfying 
$\nd^{\mathcal{T}} I=0$ where
\begin{equation}\label{trconn}
\nd^{\cT}_a
\left( \begin{array}{c}
\si\\\mu_b\\ \rho
\end{array} \right) : =
\left( \begin{array}{c}
 \nabla_a \si-\mu_a \\
 \nabla_a \mu_b+ \bg_{ab} \rho +\Rho_{ab}\si \\
 \nabla_a \rho - \Rho_{ab}\mu^b  \end{array} \right) .
\end{equation}

The connection $\nd^{\mathcal{T}}$ constructed here (following
\cite{BEG}) is the normal conformal tractor connection. We will often
write simply $\nabla$ for this when the meaning is clear by
context. This is convenient since we will couple the tractor
connection to the Levi-Civita connection.

Let us
write $J^k \ce[1]$ for the bundle of k-jets of germs of sections of
$\ce[1]$.  Considering, at each point of the manifold, sections which
vanish to first order at the given point point reveals a
canonical sequence,
$$
0\to S^2 T^*M\otimes \ce[1]\to J^2 \ce[1]\to J^1\ce[1]\to 0~.
$$ This is the jet exact sequence at 2-jets. Via the conformal metric
$\bg$, the bundle of symmetric covariant 2-tensors $S^2 T^*M$
decomposes directly into the trace-free part, which we will denote
$S^2_0 T^*M$, and a pure trace part isomorphic to $\ce [-2]$,  hence
$S^2 T^*M \otimes \ce[1]=(S^2_0 T^*M \otimes \ce[1]) \oplus \ce [-1]$.  The
{\em standard tractor bundle} $\mathcal{T}$ may defined as the quotient of
$J^2 \ce [1]$ by the image of $S^2_0 T^*M \otimes \ce[1]$ in $J^2
\ce[1]$. By construction this is invariant, it depends only on the
conformal structure. Also by construction, it is an extension of the
1-jet bundle
\begin{equation}\label{Xdef}
0\to \ce[-1]\stackrel{X}{\to} \mathcal{T} \to J^1 \ce[1]\to 0.
\end{equation}
The canonical homomorphism $X$ here will be viewed as a section of
$\mathcal{T}[1]=\mathcal{T}\otimes \ce[1]$ and, with the jet exact
sequence at 1-jets, controls the filtration structure of
$\mathcal{T}$.

Next note that there is a tautological operator $D:\ce[1]\to
\mathcal{T}$ which is simply the composition of the universal 2-jet
differential operator $j^2:\ce[1]\to \Gamma(J^2 \ce[1])$ 
followed by the canonical projection $J^2 \ce[1]\to \mathcal{T}$. 
On the other hand, via a choice of metric $g$, and the Levi-Civita
connection it determines, we obtain a differential operator $\ce[1]\to
\ce[1]\oplus \ce^1[1]\oplus \ce[-1]$ by $\si\mapsto ( \si, \nd_a \si,
 \frac{1}{d}(\Delta - \J) \si )$ and this obviously determines an
isomorphism
\begin{equation}\label{split}
\mathcal{T} \stackrel{g}{\cong} \ce[1]\oplus T^*M[1]\oplus
\ce[-1] ~.
\end{equation}
In the following we shall frequently use \nn{split}. Sometimes this
 will be without any explicit comment but also we may write for
 example $t\stackrel{g}{=}(\si,~\mu_a,~\rho)$, or alternatively
 $[t]_g=(\si,~\mu_a,~\rho)$, to mean $t$ is an invariant section of
 $\cT$ and $(\si,~\mu_a,~\rho) $ is its image  under the
 isomorphism \nn{split}.  Changing to a conformally related metric
 $\widehat{g}=e^{2\om}g$ ($\om$ a smooth function) gives a different
 isomorphism, which is related to the previous by the transformation
 formula
\begin{equation}\label{transf}
\widehat{(\si,\mu_b,\rho)}=(\si,\mu_b+\si\Up_b,\rho-\bg^{bc}\Up_b\mu_c-
\tfrac{1}{2}\si\bg^{bc}\Up_b\Up_c),  
\end{equation}
where $\Upsilon:=d\om$. It is straightforward to verify that the
right-hand-side of \nn{trconn} also transforms in this way and hence
$\nd^\mathcal{T}$ gives a conformally invariant connection on
$\mathcal{T}$ which we shall also denote by $\nd^\mathcal{T}$. This is
the tractor connection.  There is also a conformally invariant {\em
tractor metric} $h$ on $\mathcal{T}$ given (as a quadratic form) by
\begin{equation}\label{trmet}
(\si,~\mu,~\rho)\mapsto \bg^{-1}(\mu,\mu)+2\si \rho~.
\end{equation}
 This is
preserved by the connection and clearly has signature $(d+1,1)$. 

Let us return to our study of the equations \nn{primc} and \nn{ascE}.
First observe that, given a metric $g$,
through \nn{split} the tautological invariant operator $D$ from above
is given by the explicit formula
\begin{equation}\label{Dexpl}
 D:\ce[1]\to \mathcal{T} \quad \quad \quad \si\mapsto (\si,~\nd_a \si
,~ \frac{1}{d}(\Delta \si -\J \si)).
\end{equation} 
This is a differential splitting operator, since it is inverted by the
canonical tractor $X$: $h(X,D\si)=\si$. (To see this one may use that
in terms of the splitting \nn{split} $X=(0,~0,~1)$.)  If a standard
tractor $I$ satisfies $I=D\si$ for some $\si\in \ce[1]$ then
$\si=h(X,I)$ and we shall term $I$ a {\em scale tractor}. 
For the study of scale tractors the following result is useful. 
\begin{lemma}\label{zeropt} For $\si$ a section of $\ce[1]$ we have 
\begin{equation}\label{I2}
|D\si|^2:=h(D \si ,D\si )= \frac{2}{d} \si (\Delta^g-J^g)\si + |\nd^g
 \si |^2_{\bg},
\end{equation}
where $|\nabla \si |_{{\sbg}}^2$ means $\bg^{ab} (\nd_a \si)\nd_b\si$.
In particular,  if $\si(p)=0$, $p\in M$, then
$$
|D\si|^2(p) =|\nabla \si |_{{\sbg}}^2(p).
$$
\end{lemma}
\noindent{\bf Proof:} This follows easily from the formulae \nn{trmet}
and \nn{Dexpl}.
\quad $\Box$

Using Lemma \ref{zeropt}, we have the following.
\begin{proposition}\label{ascP}
If $I$ is a scale tractor then 
$$
|I|^2= - S(\si), 
$$
where $\si=h(X,\si)$. In particular off 
the zero set of
$\si$ we have
$$
|I|^2=-\frac{2}{d} \ul{J}^{\g} 
$$
where $\g=\si^{-2}\bg$ and $\ul{J}^{\g} $ is the $\g$ trace of $P^{\g}$.
An ASC structure is a conformal manifold $(M,[g])$ equipped with a
 scale tractor of constant length. 
\end{proposition}
\noindent{\bf Proof:}
Everything is clear except the point made in the second display.
Recall that now, in contrast to the Introduction,  
$ J^{\g}$ denotes $\bg^{ab}P_{ab}^{\g}$. So,
writing $g_o$ for the inverse to $\g$, we have
$$
\si^2 J^{\g}=\si^2\bg^{ab}P_{ab}^{\g}= g_o^{ab}P_{ab} =: \ul{J}^{\g},
$$ this is the metric $\g$ trace of the Schouten tensor $P^{\g}$. On
the other hand, away from the zero set of $\si$, we may calculate in
the scale $\si$ and we have $\nabla^{\g} \si=0$, whence $-2\si^2
J^{\g}/d$ is exactly the right hand side of \nn{I2}.  \quad $\Box$

Now collecting our observations we obtain the basic elements of the tractor picture for AE structures, as follows.
\begin{theorem}\label{key}
A directed almost Einstein structure is a conformal manifold $(M^{n+1},[g])$
 equipped with a parallel (standard) tractor $I\neq 0$. The mapping
 from non-trivial solutions of \nn{primc} to parallel tractors is by
 $\si\mapsto D\si$ with inverse $I\mapsto \si:=h(I,X)$. 
If $I\neq 0$ is parallel and  $\si:=h(I,X)$ then
the structure 
$(M,[g],\si)$ is ASC with $S([g],\si)=-|I|^2$. On the open set where $\si$ is nowhere vanishing
 $\g:=\si^{-2}\bg$ is Einstein with $\Ric^{\g}=n|I|^2\g$. 
\end{theorem}
\noindent{\bf Proof:} The first observation is immediate from the
construction in \nn{trconn} of the tractor connection as a
prolongation of the equation \nn{primc} for an almost Einstein
structure. 

Next observe that if $I\stackrel{g}{=}(\si,~\mu_a,~\rho)$ is a
parallel section for $\nd^{\mathcal{T}} $ then it follows immediately
from the formula \nn{trconn} that necessarily
\begin{equation}\label{pD}
 \big(\si,~\mu_a,~\rho\big) =
(\si,\nd_a \si , \frac{1}{d}(\Delta \si -J \si)), 
\end{equation} 
that is $I$ is a scale tractor, $I=D\si$. From the formula for the
tractor metric it follows that $\si=h(X,I)$. 

Since the tractor connection preserves the tractor metric it follows
that if $I$ is a parallel tractor then $|I|^2:=h(I,I)$ is
constant. Thus an almost Einstein structure is ASC as claimed.

For the final statement we use that $I$ parallel implies that $\si$
satisfies \nn{primc}.  On the set where $\si$ is nowhere vanishing we
may use the metric $\g=\si^{-2}\bg$. The corresponding Levi-Civita
connection annihilates $\si$ and then \nn{primc} asserts that $P^{\g}$
is trace-free.  \quad $\Box$ \\ In view of the the Theorem we shall
often use the notation $(M,[g],I)$ to denote a directed almost
Einstein manifold. In this context $I$ should be taken as parallel and
non-zero.

There is a useful immediate consequence of the Theorem, as follows.
\begin{corollary}\label{betw} On a fixed conformal structure $(M,[g])$ 
the set of directed AE structures is naturally a vector space with the origin
removed.  In particular if $I_1$ and $I_2$ are two linearly
independent directed AE structures then for each $t\in \mathbb{R}$ 
$$
I_t:=(\sin t)I_1+ (\cos t)I_2
$$ is a directed AE structure. In this case given $p\in M$ there is $ t\in
\mathbb{R}$ so that $\si_t(p):=h(X,I_t)_p=0$.
\end{corollary}
One might suspect that generically non-scalar positive AE manifolds
will have non-empty scale singularity sets. The Corollary shows that
this certainly is the case on a fixed conformal structure with two
linearly independent AE structures.

\section{Classification of the scale singularity set} \label{class}

Given a standard tractor $I$ and $\si:=h(X,I)$ let us write $S(I)$
as an alternate notation for $S(\si)$.  
As before we write 
$$\Sigma :=\{p\in
M~|~ \si(p)=0\}
$$ and term this the {\em scale singularity set} of $I$; this is the
set where $\g=\si^{-2}\bg$ is undefined.  In this section we shall
establish the following, and then complete to a proof of Theorem
\ref{maina}.
\begin{theorem}\label{classthm} Let $(M,[g],I)$ be an almost Einstein 
structure.  There are three cases:\\ $\bullet$ $|I|^2<0$, which is
equivalent to $S(I)>0$, then $\Sigma$ is empty and $(M,\si^{-2}\bg)$
is Einstein with positive scalar curvature; \\ $\bullet$  $|I|^2=0$,
which is equivalent to $S(I)=0$, then $\Sigma$ is either empty or
consists of isolated points, and $(M\setminus \Sigma,\si^{-2}\bg)$ is
Ricci-flat; \\
$\bullet$  $|I|^2> 0$, which is
equivalent to $S(I)<0$, then the scale singularity set $\Sigma$ is
either empty or else is a totally umbillic hypersurface, and
$(M\setminus \Sigma, \si^{-2} \bg)$ is Einstein of negative scalar
curvature.
\end{theorem}
\noindent The curvature statements follow from Theorem \ref{key}. Also
from there we have that an AE manifold is ASC. Thus from Theorem
\ref{mainsc} we have at once both the first result and also that if,
alternatively, $|I|^2>0$ then the singularity set is either empty or
is a hypersurface. The proof is completed via Propositions \ref{nullp}
and \ref{aEu} below.

We shall make a general observation which sheds light on the scalar
flat case.  From Theorem \ref{key}, $I$ parallel implies $ I= D \si$,
for some density $\si$ in $\ce[1]$.  An obvious question is whether,
at any point $p\in M$, we may have $j^1_p\si=0$, i.e.\ whether the
1-jet of $\si$ may vanish at $p$. Evidently this is impossible if
$|I|^2\neq 0$. We observe here (cf. \cite{GoPrague}) that, in any
case, if $ I= D \si\neq 0$ is parallel then the zeros of $j^1\si$ are
isolated. In fact we have a
slightly stronger result. As usual here we write $\si=h(X,I)$.
\begin{lemma}\label{j1facts} 
Suppose that $ I\neq 0$ is parallel and $j^1_p\si=0$. Then
there is a neighbourhood of $p$ such that, in this neighbourhood,
$\si$ is non-vanishing away from $p$.
\end{lemma}
\noindent{\bf Proof:} 
Suppose that
$ I\neq 0$ is parallel
and $j^1_p\si=0$. Since $I$ is parallel $ I= D \si$. This with 
$j^1_p\si=0$ implies that, at $p$, and in the scale $g$, we have 
$I\stackrel{g}{=}(0~,0~,\rho)$ for some density $\rho$ with $\rho(p)\neq 0$. 
Thus from \nn{trconn} (or equivalently \nn{primc}) we have $(\nd_a
\nd_b \si) (p)=-\rho(p) \bg_{ab}(p)$.  Trivialising the density
bundles using the metric $g$ the latter is equivalent to $(\nd_a \nd_b
s) (p)=- r(p) g_{ab}(p)$ where the smooth function $r$ satisfies
$r(p)\neq 0$.  (Here we use that $g=\tau^{-2} \bg$ for some
non-vanishing $\tau$ in $\ce[1]$ and $s=\tau^{-1}\si$ while $r=\tau
\rho$. Then since $\nd$ is the Levi-Civita for $g$ we have
$\nd\tau=0$.)  So, in terms of coordinates based at $p$, the first
non-vanishing term in the Taylor series for $s$ (based at $p$) is $-r
g_{ij}x^ix^j$.  \quad $\Box$\\

\medskip

Note that an ASC structure is scalar flat if and only if  
$\g $ is Ricci-flat on
 $M\setminus \Sigma$. In the following $\si:=h(X,I)$.
\begin{proposition}\label{nullp}
If $(M,[g],I)$ is an ASC structure with $j^1_p
\si=0$, at some point $p$, then $(M,[g],I)$ is scalar flat.
Conversely if $(M,[g],I)$ is a scalar flat ASC structure then, at any
 $p\in M$ with $\si (p)=0$ we have $j^1_p\si=0$.

If $(M,[g],I)$ is a scalar flat AE structure then, at any
 $p\in M$ with $\si (p)=0$ we have $j^1_p\si=0$
 and $j^2_p \si \neq
 0$. For any AE manifold the scale singularity set consists of isolated points.
\end{proposition}
\noindent{\bf Proof:} Since by definition $I=D\si$, from Lemma
\ref{zeropt} it is immediate that, at any point $p$ with $\si(p)=0$,
we have $S(\si)(p)=0$ if and only if $j^1_p\si=0$. (Alternatively this
is visible directly from the formula \ref{ascC}.)  The first two
statements follow immediately, as by definition $S(\si)$ is constant
on an ASC manifold.

Now we consider AE manifolds. These are ASC and so we have the first
 results.  Since $I$ is parallel, we have $I= D\si$.  If an AE
 manifold is scalar flat then, at a point $p$ where $\si(p)=0$, we
 have $j^1_p\si=0$ and so from \nn{pD} the tractor $I$ is of the form
 $I\stackrel{g}{=}(0,~0,~ \rho)$ at $p$.
On the other hand, since $I\neq 0$ is parallel, it
follows that $D\si=I$ is nowhere zero on $M$. Hence (since $D$ is a
second order differential operator) $j^2\si$ is non-vanishing.  In
fact, from \nn{Dexpl}, at any point $p$ where $j^1_p\si$ vanishes we
have $\rho(p)=\frac{1}{d} (\Delta \si)(p) \neq 0$. The last statement is now an immediate
consequence of Lemma \nn{j1facts}.  \quad $\Box$\\

\noindent{\bf Remark:} Note that $j^1_p\si=0$ means that when we work
in terms of a background metric $g$ we have $j^1_ps=0$ for the
function $s$ corresponding to $\si$ and so $p$ is a critical point of
$s$. In fact it is already clear from \nn{ascE} that, even for ASC
structures, if $S(g,s)=0$ then $s_p=0$ implies $p$ is a critical
point.

\subsection{Conformal hypersurfaces and the scale singularity set}\label{hyp}

Let us first recall some facts concerning  general hypersurfaces in a
 conformal manifold $(M^d,[g])$, $d\geq 3$. 
 If $\Sigma$ is a boundary component of a Riemannian (or conformal)
 manifold then, without further comment, we will assume that the
 conformal structure extends smoothly to a collar of the boundary. Our
 results will not depend on the choice of extension.  So in the
 following we suppose that $\Sigma$ is an embedded codimension 1
 submanifold of $M$.

 Let $n_a$ be a
section of $\ce_a[1]$ such that, along $\Sigma$, we have
$|n|^2_{\sbg}:=\bg^{ab}n_a n_b=1$. Note that the latter is a conformally
invariant condition since $\bg^{-1}$ has conformal weight $-2$.  Now in the
scale $g$, the mean curvature of $\Sigma$ is given by
$$
H^g=\frac{1}{n-1}\big(\nd_a n^a -n^an^b \nd_an_b  \big),
$$ 
as a conformal $-1$-density. This is independent of how $n_a$ is
extended off $\Sigma$.  Now under a conformal rescaling, $g\mapsto
\widehat{g}=e^{2\om} g$, $H$ transforms to $
H^{\widehat{g}}=H^g+n^a\Upsilon_a$. Thus we obtain a conformally invariant
section $N$ of $\cT|_\Sigma$
$$
N\stackrel{g}{=}\left(\begin{array}{c}0\\
n_a\\
-H^g\end{array}\right),
$$ and from \nn{trmet} $h(N,N)=1$ along $\Sigma$. Obviously $N$ is
independent of any choices in the extension of $n_a$ off
$\Sigma$. This is the {\em normal tractor} of \cite{BEG} and may be
viewed as a tractor bundle analogue of the unit conormal field from
the theory of Riemannian hypersurfaces.

Recall that a point $p$ in a hypersurface is an umbillic point if 
at that point the second fundamental form is trace free, this is a
conformally invariant condition.  A hypersurface is totally
umbillic if this holds at all points.  
Differentiating $N$ tangentially along $\Sigma$ using $\nd^\cT$, directly 
from \nn{trconn} we obtain 
the following
result. 
\begin{lemma}\label{umbillic}
If the normal tractor $N$ is parallel, with respect to $\nabla^{\cT}$, along a hypersurface $\Sigma$
then  the hypersurface $\Sigma$ is totally umbillic.
\end{lemma}
\noindent In fact constancy of $N$ along a hypersurface is equivalent to total
umbillicity. This is (Proposition 2.9) from \cite{BEG}.

Now let us return to the study of ASC and AE structures. First we see
that the normal tractor is linked, in an essential way, to the
ambient geometry off the hypersurface.
\begin{proposition}\label{ascN}
Let $(M^d,[g],I)$ be a scalar negative ASC
structure with scale singularity set $\Sigma\neq \emptyset$ and
$|I|^2=1$. Then, with $N$ denoting the normal tractor for $\Sigma$, we
have $N=I|_\Sigma$.
\end{proposition}
\noindent{\bf Proof:}
As usual let us write $\si:=h(X,I)$.  By definition 
$$
I=D\si\stackrel{g}{=} 
\left(\begin{array}{c} 
\si \cr  \nabla_a \si\cr \frac{1}{d}(\Delta \si - J \si) 
\end{array}\right) ~.
$$ 
 Let us write $n_a := \nabla_a \si$.  Along
$\Sigma$ we have $\si=0$, therefore 
 $$
I|_\Sigma  \stackrel{g}{=} \left(\begin{array}{c}
 0 \cr  n_a \cr
 \frac{1}{d}\Delta \si \end{array}\right)  ~,
$$ and from Lemma \ref{zeropt} $|n|^2_{\sbg}=1$, since $|I|^2=1$. So
$n_a|_\Sigma$ is a conformal weight 1 conormal field for $\Sigma$.

Next we calculate the mean curvature $H$ in terms of $\si$.  Recall $
(d-1)H=\nd^an_a-n^an^b\nd_b n_a , $ on $\Sigma$. 
 We calculate the right hand side in
a neighbourhood of $\Sigma$. Since $n_a=\nabla_a\si$, we have
$\nd^a n_a=-\Delta \si$. On the other hand
$$
n^an^b\nd_b n_a=\frac{1}{2}n^b\nd_b
(n^an_a)=\frac{1}{2} n^b\nd_b(1-\frac{2}{d}\si\Delta \si +\frac{2}{d}J\si^2),
$$
where we used that $|D \si|^2=1$.
Now along $\Sigma $ we have $1=n^an_a=n^a\nd_a \si$, and so there this simplifies to 
$$
n^an^b\nd_b n_a=-\frac{1}{d}\Delta \si.
$$
Putting these results together, we have
$$
(d-1)H= \frac{1}{d}(1-d)\Delta \si\quad \Rightarrow \quad H=-\frac{1}{d}\Delta \si~.
$$

Thus 
$$
I|_\Sigma  \stackrel{g}{=} \left(\begin{array}{c}
0 \cr  n_a \cr
 -H \end{array}\right) ~,
$$ as claimed. 
 \quad $\Box$

A consequence for AE structures follows easily.
\newcommand{\II}{I \hspace*{-3pt} I}
\begin{proposition}\label{aEu}
Let $(M^d,[g],I)$ be a scalar negative almost Einstein structure with
scale singularity set $\Sigma\neq \emptyset$ and $|I|^2=1$. Then
$\Sigma $ is a totally umbillic hypersurface with $I|_\Sigma=N$, the
normal tractor for $\Sigma$.
\end{proposition}
\noindent{\bf Proof:} Since an AE structure $(M,[g],I)$
is ASC it follows from Proposition \ref{ascN}
above that along the singularity hypersurface $I$
agrees with the normal tractor $N$.    On
the other hand, since $I$ is parallel everywhere, it follows that $N$
is parallel along $\Sigma$ and so, from Lemma \ref{umbillic}, $\Sigma$
is totally umbillic.  \quad $\Box$

Proposition 2.8 of LeBrun's \cite{LeBrunH} also gives a proof that the
conformal infinity of a PE metric is totally umbillic. 

\medskip

\noindent{\bf Proof of Theorem \ref{classthm}:} The remaining point is
to show that if $(M,[g],I)$ is AE with $|I|^2>0$ and a singularity
hypersurface $\Sigma$, then this is totally umbillic.  This is
immediate from the previous Proposition as multiplying $I$ with a
positive constant yields a 
yields a parallel tractor with the same singularity set.
\quad $\Box$

\medskip

Most of Theorem \ref{maina} is simply repackaging of the tractor based
statements in Theorem \ref{classthm} above. To complete the proof of
the former we simply need to describe PE manifolds in the same
language, and this is our final aim for this section.
\begin{proposition}\label{peae}
Suppose that  $M$  is a compact manifold with boundary
$\Sigma$, and $(M,[g],I)$ is an almost Einstein structure with $|I|^2=1$, and
such that the scale singularity set is $\Sigma$. Then $(M,[g],I)$ is a
Poincar\'e-Einstein manifold with the interior metric $\g=\si^{-2}\bg$, where 
$\si:= h(X,I)$.
Conversely
Poincar\'e-Einstein manifolds are scalar negative almost Einstein
structures. 
\end{proposition}

\noindent{\bf Proof:} Suppose that $(M,[g],I)$ is an AE structure as
described.  Since AE manifolds are ASC, with the parallel tractor $I$
giving the scale tractor of the ASC structure, it follows from Theorem
\ref{mainsc} that $(M,[g],\si)$ is AH. But $I$ parallel means that
$\g=\si^{-2}\bg$ is Einstein on $M\setminus \S$, and there $|I|^2=1$
is equivalent to $\Ric (g^o)=-n g^o$.
The converse direction is also straightforward, or see \cite{GoMGC}. 
\quad $\Box$

\section{Conformal geometry of $\Sigma$ versus conformal geometry of $M$}
\label{MvS}

Here for almost Einstein manifolds we shall derive basic equations
satisfied by the conformal curvatures.
In particular
 for Poincar\'e-Einstein manifolds, and more generally for scalar
negative almost Einstein manifolds, we shall study the relationship between
the conformal geometry of $M$ and the intrinsic conformal geometry of
the scale singularity set $\Sigma$.  Since $\Sigma$ is a hypersurface,
a first step is to understand the conformal structure induced on an
arbitrary hypersurface in a conformal manifold and in particular the
relationship between the intrinsic conformal tractor bundle of $\S$ and
the ambient tractor bundle of $M$. This is the subject of Section
\ref{hyp2}. On the other hand we have already observed that on scalar
negative AE manifolds the singularity set is umbillic. So the main aim
of this section is to deepen this picture. We shall see that the along
the singularity hypersurface the intrinsic tractor connection
necessarily agrees with an obvious restriction of the ambient tractor
connection.  This has immediate consequences for the relationship
between the intrinsic and ambient conformal curvature quantities, but
we are able to also show that there is an even stronger compatibility between 
the conformal curvatures of $(\S,[g_\S])$ and those of $(M,[g])$.
Finally we shall derive equations on the latter that partly establish  a
Dirichlet type problem based directly on the  conformal curvature quantities.

\subsection{Conformal hypersurfaces}\label{hyp2}

Here we revisit (cf.\ Section \ref{hyp}) the study of a general
 hypersurface $\Sigma$ in a conformal manifold $(M^d,[g])$, $d\geq
 3$. This time our aim is to see, in this general setting, how the
 conformal structure of the hypersurface is linked that of the ambient
 space.

 With respect to the embedding map, each metric $g$ from the conformal
 class on $M$ pulls back to a metric $g_\Sigma $ on $\Sigma$. Thus the
 ambient conformal structure of $M$ induces a conformal structure
 $[g_\Sigma]$ on $\Sigma^n$ ($n+1=d$); we shall refer to this as the
 {\em intrinsic} conformal structure of $\Sigma$. Given the
 relationship of the intrinsic and ambient conformal structures it
 follows easily that the intrinsic conformal density bundle of weight
 $w$, $\ce^\S[w]$ is canonically isomorphic to $\ce[w]|_\S$ and we
 shall no longer distinguish these.  It is also clear that since
 $g_\S$ is determined by $g$ the trivialisations they induce on,
 respectively, $\ce^\S[w]$ and $\ce[w]$ are consistent. In particular
 the Levi-Civita connection on $\ce^\S[w]$ from $g^\S$ agrees with
 the restriction of the connection on $\ce[w]$ arising from the
 trivialisation due to $g$.

If $n\geq 3$ then $(\Sigma, [g_\Sigma])$ has an intrinsic tractor
bundle $\cT_\S$. We want to relate this to $\cT$ along $\Sigma$. Note
that $\cT_\S$ has a canonical rank $n+2$ subbundle, viz.\ $N^\perp$ the
orthogonal complement (with respect to $h$) of the normal tractor $N$. 
As noted in \cite{BrGonon}, there is a canonical isomorphism 
\begin{equation}\label{trisom}
N^\perp \stackrel{\cong}{\longrightarrow} \ct_\S ~.
\end{equation} 
To see this let $n_a$ denote a weight 1 conormal field along $\S$. 
There is a canonical inclusion of $T\S$ in $TM|_\Sigma$ and 
we
identify $T^*\S$ with the annihilator subbundle
in  $T^*M|_\S$ of $n^a$. These identifications do not
require choosing a metric from the conformal class. Now calculating
in a scale $g$ on $M$, $\cT$ and hence also $N^\perp$, decomposes into
a triple via \nn{split}. Then the mapping of the isomorphism 
is (cf.\ \cite{Grant})
\begin{equation}\label{xmap}
[N^\perp]_g\ni \left( \begin{array}{c}
\si\\
\mu_b\\ \rho
\end{array} \right) \mapsto 
\left( \begin{array}{c}
\si\\
\mu_b  -H n_b \si
\\ \rho +\frac{1}{2}H^2 \si
\end{array} \right) \in [\cT_\S]_{g_\S}
\end{equation}
 where, as usual, $H$ denotes the mean curvature of $\S$ in the
scale $g$ and $g_\S$ is the pullback of $g$ to $\S$.  Since
$(\si,\mu_b,\rho )$ is a section of $[N^\perp]_g$ we have $n^a \mu_a=
H\si$. Using this one easily verifies that the mapping is conformally
invariant: If we transform to $\widehat{g}=e^{2\om}g$, $\om\in \ce$,
 then $(\si,\mu_b,\rho )$
transforms according to \nn{transf}. Using that
$\widehat{H}=H+n^a\Upsilon_a$ one calculates that the image of
$(\si,\mu_b,\rho )$ (under the map displayed) transforms by the
intrinsic version of \nn{transf}, that is by \nn{transf} except where
$\Upsilon_a$ is replaced by $\Upsilon^\S_a=\Upsilon_a -
n_an^b\Upsilon_b$ (which on $\S$ agrees with $d^\S \om$, the intrinsic
exterior derivative of $\om$). This signals that the explicit map
displayed descends to a conformally invariant map \nn{trisom}.

So far we understand the tractor bundle on $\S$ for $n\geq 3$.  In the
case of $n=2$, $\S$ does not in general have a preferred intrinsic conformal
tractor connection. There is much to be said in this case but for our current
purposes it will be most economical to proceed as follows.
We shall {\em define} $\ct_\Sigma$ to be the orthogonal
complement of $N$ in $\ct|_\Sigma$  and in any dimension $d\geq 3$
let us write $\Proj_{\Sigma}:\ct|_\Sigma\to \ct_\Sigma$ for the
orthogonal projection afforded by $N$. 
Then for $d=3$, equivalently $n=2$, we define the tractor connection on
$\Sigma$ to be the orthogonal projection of the ambient tractor connection.
That is, working locally, for $v\in \Gamma(T\Sigma)$ and $T\in
\ct_\Sigma=N^\perp$ we extend these smoothly to $v\in \Gamma(TM)$ and $T\in
\ct$. Then we define $\nabla^{\ct_\Sigma}_v T
:=\Proj_\Sigma(\nabla^\ct_v T)$ along $\Sigma$. It is verified by
standard arguments that this is independent of the extension choices
and defines a connection on $\ct_\Sigma$.

\medskip

Finally we observe a useful alternative approach to the arguments
above via a result that, for other purposes, we will call on later.
\begin{proposition}\label{H0}
Let $\S$ be an orientable hypersurface in an orientable conformal
manifold $(M,[g])$. In a neighbourhood of $\Sigma$ there is a metric
$\widehat{g}$ in the conformal class so that $\S$ is minimal, i.e.\ $
H^{\widehat{g}}=0.  $
\end{proposition}
\noindent{\bf Proof:} For simplicity let calculate in the metric $g$
and write $H^g$ to be the mean curvature of $\S$ as a {\em function}
along $\S$. Take any smooth extension of this to a function on $M$.
By a standard argument one can show that in a neighbourhood of
$\Sigma$ there is a {\em normal} defining function $s$ for $\Sigma$,
that is \ $\S$ is the zero set of $s$, and along $\S$ the 1-form $ds$
satisfies $|ds|^2_g=1$.  
 Then $n^a:=g^{ab}\nabla_b s$ is a unit normal vector field along $\S$.
Recall the conformal transformation of the mean curvature: If
$\widehat{g}=e^{2\om}g$, for some $\om\in \ce$ then $e^{\om}
H^{\widehat{g}}=H^g+n^a\Upsilon_a =H^g+n^a\nabla_a \om $.  Thus if
we take $\om:= - s H^g$ then $H^{\widehat{g}}=0$.  \quad $\Box$

 Dropping the `hat' on $\widehat{g}$, we see that with such
$g$ (satisfying $H^g=0$) the map \nn{xmap} simplifies significantly in this normalisation;
the splittings of $N^\perp$ and $\cT_\S$ then agree in the ``obvious
way''. This is consistent with conformal transformation: The condition
$H=0$ does not fix the representative metric $g$, even along $\S$. For
example at the 1-jet level the remaining freedom along $M$ is to
conformally rescale by $g\mapsto e^{2\om}g$ where $n^a\nabla_a
\om=0$. This is exactly as required to preserve the agreement of the
splittings of $N^\perp$ and $\cT_\S$. In fact this was the point of
view taken in \cite{BrGonon}. From there one easily recovers the
formula \nn{xmap}.

Finally we note here that the rescaling
involved in the proof of the proposition above is global and
especially natural in the case of directed ASC and AE structures.
\begin{corollary}\label{minc}
Let $(M,[g],I)$ be a directed scalar negative ASC manifold with a
scale singularity set. Then there is a metric $\widehat{g}\in [g]$
with respect to which $\S$ is a minimal hypersurface. In particular if
$(M,[g],I)$ is a directed AE manifold then $\S$ is totally geodesic
with respect to $\widehat{g}$.
\end{corollary}
\noindent {\bf Proof:} Suppose that $\si$ is the conformal weight 1
density defining a $S(\si)=-1$ ASC structure with a non-trivial scale
singularity hypersurface $\S$. Write $H^g$ (now as a $-1$ density) for
the mean curvature of $\S$ with respect to an arbitrary background
metric $g$ and extend this smoothly to $M$. Then $\S$ has mean
curvature zero with respect to the metric $\widehat{g}=e^{2\om}g$
where $\om:= -H^g \si$. For the last statement we recall that if $\si$
satisfies \nn{primc} then $\S$ is totally umbillic and this is a
conformally invariant condition. \quad $\Box$

\subsection{Tractor curvature}\label{moretr}

We digress briefly to recall some further background. In this section
we work on an arbitrary conformal manifold $(M^d,[g])$.  It will be
convenient to introduce the alternative notation $\ce^A$ for
the tractor bundle $\cT$ and its space of smooth sections. Here the
index indicates an abstract index in the sense of Penrose and so we
may write, for example, $V^A\in \ce^A$ to indicate a section of the
standard tractor bundle. Using the abstract index notation the tractor
metric is denoted $h_{AB}$ with inverse $h^{BC}$. These will be used
to lower and raise indices in the usual way.

In computations, it is often
useful to introduce the `projectors' from $\ce^A$ to the components
$\ce[1]$, $\ce_a[1]$ and $\ce[-1]$ which are determined by a choice of
scale.  They are respectively denoted by $X_A\in\ce_A[1]$,
$Z_{Aa}\in\ce_{Aa}[1]$ and $Y_A\in\ce_A[-1]$, where
$\ce_{Aa}[w]=\ce_A\otimes\ce_a\otimes\ce[w]$, etc.  Using the metrics
$h_{AB}$ and $\bg_{ab}$ to raise indices, we define $X^A, Z^{Aa},
Y^A$. Then we immediately see that
$$
Y_AX^A=1,\ \ Z_{Ab}Z^A{}_c=\bg_{bc}
$$
and that all other quadratic combinations that contract the tractor
index vanish.


Given a choice of conformal scale we have the corresponding
Levi-Civita connection on tensor and density bundles and we can use
the coupled Levi-Civita tractor connection to act on sections of the
tensor product of a tensor bundle with a tractor bundle and so
forth. This operation is defined via the Leibniz rule in the usual way.  In
particular we have
\begin{equation}\label{connids}
\nd_aX_A=Z_{Aa},\quad
\nd_aZ_{Ab}=-\Rho_{ab}X_A-Y_A\bg_{ab}, 
\quad \nd_aY_A=\Rho_{ab}Z_A{}^b.
\end{equation}

The curvature $ \Omega$ of the tractor connection 
is defined by 
\begin{equation}\label{curvature}
[\nd_a,\nd_b] V^C= \Omega_{ab}{}^C{}_EV^E 
\end{equation}
for $ V^C\in\ce^C$.  Using
\eqref{connids} and the usual formulae for the curvature of the
Levi-Civita connection we calculate (cf. \cite{GoPetCMP})
\begin{equation}\label{tractcurv}
\Omega_{abCE}= Z_C{}^cZ_E{}^e C_{abce}-X_{C}Z_{E}{}^e A_{eab}+ X_{E}Z_{C}{}^e A_{eab} 
\end{equation}
where 
\begin{equation}
A_{abc}:=2\nabla_{[b}\Rho_{c]a}  \label{cot} 
\end{equation}
is the {\em Cotton tensor}.

Next we note that there is a conformally invariant differential
operator between weighted tractor bundles
$$
\bD_A\colon\ce_{B \cdots E}[w]\to\ce_{AB\cdots E}[w-1],
$$ given a choice of conformal scale $g$, the {\em tractor-$D$
operator} by
\begin{equation}\label{Dform}
\bD_A V:=(d+2w-2)w Y_A V+ (d+2w-2)Z_{Aa}\nabla^a V +X_A (\Delta - w J) V. 
\end{equation} 
This is the (Thomas) tractor-D operator as recovered in \cite{BEG};
see \cite{GoSrni99,GoAdv} for an invariant derivation.  The conformal
operator $D$ from Section \ref{al} is simply $\frac{1}{d}$ times $\bD$
applied to $\ce[1]$. (It is convenient to retain the two notations,
rather than carry the factor $1/d$ into many calculations.) Using $\bD$
we obtain (following \cite{GoSrni99,GoAdv}) a conformally invariant
curvature quantity as follows
\begin{equation}\label{Wdef}
W_{BC}{}^E{}_F:=
\frac{3}{d-2}\bD^AX_{[A} \Omega_{BC]}{}^E{}_F ,
\end{equation}
where $\Omega_{BC}{}^E{}_F:= Z_A{}^aZ_B{}^b \Om_{bc}{}^E{}_F$.
In a choice of conformal scale, 
 $W_{ABCE}$ is given by
\begin{equation}\label{Wform}
\begin{array}{l}
(d-4)\left( Z_A{}^aZ_B{}^bZ_C{}^cZ_E{}^e C_{abce}
-2 Z_A{}^aZ_B{}^bX_{[C}Z_{E]}{}^e A_{eab}\right. \\ 
\left.-2 X_{[A}Z_{B]}{}^b Z_C{}^cZ_E{}^e A_{bce} \right)
+ 4 X_{[A}Z_{B]}{}^b X_{[C} Z_{E]}{}^e B_{eb},
\end{array}
\end{equation}
where 
\begin{equation}\label{Bdef}
B_{ab}:=\nabla^c
A_{acb}+\Rho^{dc}C_{dacb}.
\end{equation}
is known as the {\em Bach tensor} or the Bach curvature. From the formula \nn{Wform} it
is clear that $W_{ABCD}$ has Weyl tensor type symmetries.  It is shown
in \cite{CapGoamb} and \cite{GoPetCMP} that the tractor field
$W_{ABCD}$ has an important relationship to the ambient metric of
Fefferman and Graham. See also Section \ref{extS} below.

For later use we recall here some standard identities which arise from
the Bianchi identity $\nabla_{a_1}R_{a_2a_3 de}=0$, where sequentially
labelled indices are skewed over:
\begin{equation} \label{WBianchiid} 
 \nd_{a_1}C_{a_2a_3 cd} =g_{ca_1}A_{da_2a_3}-g_{da_1}A_{ca_2a_3};
\end{equation} 
\begin{equation} (n-3)A_{abc}=\nabla^d
C_{dabc}; \label{bi1} 
\end{equation} 
\begin{equation} 
\nabla^a\Rho_{ab}=\nabla_b \J ; \label{bi2} 
\end{equation}
\begin{equation} \nabla^a A_{abc}=0.\label{bi3} 
\end{equation}

\subsection{Further geometry of the singularity set} \label{depometry}
\
We are now set to return to the almost Einstein setting.
Via the projectors, a  general tractor $I^A\in \ce^A$ expands to  
$$
I^E =  Y^E\si + Z^{Ed}\mu_d+X^E\rho,
$$
where, for example, $\si=X_A I^A$. 
Hence 
$$
\Omega_{abCE} I^E= \si Z_{C}{}^c A_{cab}  +Z_C{}^c \mu^d C_{abcd}- X_{C}\mu^d A_{dab} .
$$ Now assume that $I^A\neq 0$ is parallel (of any length). As a point
of notation: in this case we shall write $I^E = Y^E\si +
Z^{Ed}n_d+X^E\rho $. That is $n_d=\nd_d\si$.  Then the left-hand-side of the last display
vanishes, whence the coefficients of $Z_C{}^c$ and $X_C$ must vanish,
i.e.,
$$
 \si A_{cab} + n^d C_{abcd} =0  \quad \mbox{and}\quad n^d A_{dab} =0.
$$ Away from the zero set of $\si$, we have that
$\si^{-1}n^d=\si^{-1}\nabla^d \si$ is a gradient and the 
first equation of the display is the condition that the metric is
conformal to a Cotton metric  (cf.\ e.g. \cite{Besse,GoNur,KNT}). On the
other hand at a point $p$ where $\si(p)=0$ the same equation shows that
\begin{equation}\label{Cn0}
C_{abcd}\nabla^d \si=C_{abcd}n^d =0  \quad \mbox{at} \quad p~.
\end{equation}

 Once again using the formulae \nn{connids} for the
tractor connection we obtain
\begin{equation}\label{divtrc}
\nd^a\Om_{acDE}= 
(d-4)Z_D{}^dZ_E{}^e A_{cde}- X_{D}Z_{E}{}^e B_{ec} + X_{E}Z_{D}{}^e B_{ec},
\end{equation}
where $B_{ab}:=\nabla^c A_{acb}+\Rho^{dc}C_{dacb}$ is the Bach tensor.
This too is annihilated by contraction with the parallel tractor $I^E$
and so we obtain
$$
(d-4)Z_D{}^d n^e A_{cde}- X_{D} n^e B_{ec}
+ 
\si Z_{D}{}^d B_{dc} =0 .
$$
From the coefficient of $Z_{D}{}^d $ we have 
$$
\si B_{dc}+(d-4) n^e A_{cde} =0 ~.
$$ In dimension four $B_{dc}$ is conformally invariant and  this
recovers the well known result that, in this dimension, it vanishes on
the conformally Einstein part of $M$. But then by continuity it
follows that the Bach tensor vanishes everywhere on $M$.  In other
dimensions  the last display shows that
$n^eA_{cde}=0 $ at any zeros of $\si$. This with \nn{Cn0} gives the
first part of the following.
\begin{proposition} \label{aEsum}
Consider an almost Einstein manifold $(M,[g],I)$
and let $\si:= I_AX^A$. We have 
$$
\begin{aligned}
 \si A_{cab} + n^d C_{abcd} =0, \quad \Rightarrow  \quad n^c A_{cab} =0,& \mbox{ and} \\
\si B_{ac}+(d-4) n^e A_{cae} =0 \quad \Rightarrow  \quad n^aB_{ab}=0,&
\end{aligned}
$$
everywhere on $M$. 
Hence for any point $p$ with $\si(p)=0$ we have   
$$
n^a C_{abcd}=0\quad \mbox{at} \quad p.
$$
In dimension $d=4$ we have $C_{abcd}(p)=0$, while $B_{ab}=0$ on $M$.
In dimensions $d\neq 4$ we have:
$$
n^a\Omega_{abCD}=0 \quad \mbox{at} \quad  p.
$$

In any dimension, if $j^1_p\si=0$ then 
$$
C_{abcd}=0=A_{bcd} \Leftrightarrow \Omega_{abCD}=0\quad \mbox{at} \quad p, \quad \mbox{and} \quad
 W_{ABCD}(p)=0.
$$
\end{proposition}
\noindent{\bf Proof:} The displayed implications follow by contracting
$n^a$ into the equations and using the symmetries of $A$ and $C$.  In
dimension 4 $n^a C_{abcd}=0$ at $p$ implies $C_{abcd}(p)=0$. When $I$
is not null (and so $n^a(p)\neq 0$), this uses the fact that we are in
Riemannian signature and is an immediate consequence of the well known
dimension 4 identity $4C_{abcd}C^{ebcd}=\delta^e_a|C|^2$. It remains
to establish the final claims.  If at some point $p$ we have
$j^1_p\si=0$, then, at $p$ we have $I^A=\rho X^A$ with $\rho(p)\neq
0$. So from $ I^E\nabla_a\Omega_{bcDE}=0 $ it follows that $X^E
\nabla_a\Omega_{bcDE}=0$ at $p$.  But, $\nd_a X^E=Z^E{}_a$ and from
\nn{tractcurv} we have $X^E\Om_{bcDE}=0$ everywhere. So $ Z_D{}^d
C_{bcda}-X_{D} A_{abc}= Z^E{}_a\Om_{bcDE}=0$ at $p$. This immediately
yields the result of the last display. But we also have that
$X^AW_{ABCD}$ and $I^A W_{ABCD}$ are zero everywhere, so an easy
variation of the last argument also shows that $W_{ABCD}$ vanishes at
$p$.  \ \quad $\Box$

In the case that $I$ is null $\g=\si^{-2}\bg$ is Ricci flat on $M\setminus \Sigma$. So, from the last part of the Proposition, it follows that $\g$ is
asymptotically flat (locally) as we approach any points of $\S$.
Following \cite{GoNur} let us say a conformal manifold of dimension
$d\geq 4$ is {\em weakly generic} at $p\in M$ if the only solution at
$p$ to $C_{abcd}v^d=0$ is $v_p^d=0$; then say that $(M,[g])$ is {\em
weakly generic} if this holds at all points of $M$. From the
Proposition above and Corollary \ref{betw} we see that the rank of the
Weyl tensor obstructs certain AE structures. Summarising we have the
following.
\begin{corollary}\label{weakgen} Let $(M^d,[g],I)$ be an almost Einstein 
structure with $S(I)=-|I|^2=0$. Then $(M,\g)$ is  asymptotically locally
Euclidean as we approach any point $p$ with $\si(p)=0$.
If $(M^{d\geq 4},[g],I)$ is an AE structure with scale singularity set
$\S\neq \emptyset$ then $(M,[g])$ is not weakly generic. If $(M,[g])$
admits any two linearly independent AE structures then it is nowhere
weakly generic.
\end{corollary}
$S(I)=0$ AE structures were studied via a different approach in
\cite{KR}; as well as some of the results mentioned here they show
that if $M$ is closed, or $(M\setminus \S,\g)$ is complete, then
$\S\neq \emptyset$ implies that $(M,g)$ is conformally diffeomorphic
to the standard sphere. They also discuss the asymptotic flatness in
preferred coordinates based at $p$.

\medskip

Now we specialise to the case of a scalar negative almost Einstein
manifold $(M,[g], I)$, with a non-empty scale singularity set
$\Sigma$.  We may suppose, without loss of generality, that $|I|^2=1$.
From Corollary \ref{minc} we may also assume that $g$ is a metric in
the conformal class so that $H^g=0$, where $H^g$ is the mean curvature of the hypersurface $\S$. 

As usual we identify $T\Sigma$ with its image in
$TM|_\Sigma$ under the obvious inclusion and $T^*\Sigma$ with the
orthogonal complement of $n_a$.  In our calculations here we will
reserve the abstract indices $i,j,k,l$ for $T\Sigma\subset TM|_\Sigma$
and its dual.  For
example $R_{ijcd}$ means the restriction of the Riemannian curvature $R_{abcd}=R^g_{abcd}$ to tangential
(to $\Sigma$) directions in the first two slots.
Now, calculating in the metric $g$,
recall that the Riemannian curvature $R_{abcd}$ decomposes
into the totally trace-free Weyl curvature $C_{abcd}$ and a remaining
part described by the Schouten tensor $P_{ab}$, according to
\nn{Rsplit}.
It follows that along  $\Sigma$
$$
R_{ijkl}=C_{ijkl}+2\bg^\Sigma_{k[i}\Rho_{j]l}+2\bg^\Sigma_{l[j}\Rho_{i]k},
$$ where we have used that the intrinsic conformal metric on $\Sigma$
is just the restriction of the ambient conformal metric. The
Levi-Civita connection $\nabla$ on $(M,g)$ induces a connection on
$T\Sigma$ (this is by differentiating tangentially followed by
orthogonal projection into $\Gamma (T\S)$). It is easily verified that
the induced connection is torsion free and on the other hand, since
$\Sigma$ is totally geodesic for $g$, it follows that the induced
connection preserves the induced metric $g_\S$.  Thus we find the
standard result that for totally geodesic hypersurfaces the induced
parallel transport agrees with the intrinsic parallel transport. It
follows immediately that $R_{ijkl}=R^{\Sigma}_{ijkl}$, where by
$R^\Sigma$ we mean the intrinsic Riemannian curvature of
$(\Sigma,g|_\S)$.  But since $n^aC_{abcd}=0$ we have that $
C_{ijkl}|_\Sigma $ is completely trace-free with respect to
$\bg^\Sigma$ and so has Weyl-tensor type symmetries, as a tensor on
$\Sigma$. It follows easily that, for $d\geq 4$, 
the right-hand-side of the last display necessarily gives the
canonical decomposition of $R^\Sigma_{ijkl}$ into its Weyl and
Schouten parts. On the other hand using again \nn{Rsplit}, but now applied to 
$R_{ijcd}$  we see that $P_{ib}n^b=0$. 
That is along $\Sigma$
\begin{equation}\label{rform}
C^\Sigma_{ijkl}=C_{ijkl}, \quad P^\Sigma_{ij}= P_{ij}\quad \mbox{and} \quad
P_{ib}n^b=0
\end{equation}
Note that since the Weyl curvature of any 3 manifold is identically
zero, in the case of dimension $d=4$ we have $C^\Sigma\equiv 0$. Thus in this
dimension the display is consistent with Proposition \ref{aEsum} where we
observed that $C|_\Sigma=0$.

As discussed in Section \ref{hyp2}, $\cT_\Sigma$ may be
identified with $N^\perp$ (i.e. the orthogonal complement of the
normal tractor) in $\cT|_\Sigma$ and we shall continue to make this
identification. Since $N$ is parallel along $\S$, $\nabla^\cT$
preserves this subbundle.  Now recall that the conformal density
bundles on $\Sigma$ are just the restrictions of their ambient
counterparts: $\ce^\Sigma[w]=\ce[w]|_\Sigma$. When we work with the
metric $g$, which has $H^g=0$, then the splittings of the tractor
bundles also coincide in the obvious way (see the Remark concluding
Section \ref{hyp2}), and in particular (via the intrinsic version of
\nn{split}) $\cT_\Sigma$ decomposes to $\ce^\Sigma[1]\oplus
\ce^\Sigma_i[1]\oplus \ce^\Sigma[-1]$ where the weight one 1-forms on
$\Sigma$, $ \ce^\Sigma_i[1]$ may be identified with $n^\perp$ in
$\ce_a[1]|_\Sigma$.  It follows from these observations, the explicit
formula \nn{trconn} expressed with respect to the metric $g$, and the
second result in the display \nn{rform}, that the tractor parallel
transport on $\Sigma^{n\geq 3}$ is just the restriction of the
ambient.  Although we used special scales for the argument it suffices
to use any metric from the conformal class to verify the agreement
since the connections are conformally invariant.  Let us summarise the
consequences.
\begin{theorem}\label{umbiltr}
Let $(M^{d\geq 3},[g],I)$ be a scalar negative almost Einstein
structure with a non-empty scale singularity hypersurface $\S$. The
 tractor connection of $(M,[g])$ preserves the intrinsic tractor
bundle of $\Sigma$, where the latter is viewed as a subbundle of the ambient
tractors: $\cT_\Sigma\subset \cT$.  Furthermore the intrinsic tractor parallel
transport of $\nabla^{\ct_\Sigma}$ coincides with the restriction of
the parallel transport of $\nabla^\cT$. 

We have 
$$
\Omega(u,v)=\Omega^\Sigma(u,v) \quad \mbox{ along } \Sigma
$$
where $u,v\in \Gamma(T\Sigma)$. In dimensions $d\neq 4$ we have the stronger result
$$
\Omega(\cdot,\cdot)=\Omega^\Sigma(\cdot,\cdot)  \quad \mbox{ along } \Sigma ,
$$ where here, by trivial extension, we view $\Omega^\Sigma$ as a section of 
$\Lambda^2 T^*M\otimes \End \cT$.
\end{theorem}
\noindent{\bf Proof:} In the case of $d=3$ the agreement of the
parallel transport is immediate from the definition of the tractor
connection $\nabla^{\ct_\Sigma}$ and that the normal tractor $N^A$ is
parallel along $\Sigma$. In the remaining dimensions this was
established immediately above. From this, and the fact that on
$\Sigma$ we have $\Omega(u,v)N=0$, it follows at once that
$\Omega(u,v)=\Omega^\Sigma(u,v)$ along $\Sigma$, as claimed.  For
dimensions $d\neq 4$ we have from Proposition \ref{aEsum} above that
$\Omega(n,\cdot)=0$, whence the final claim. \quad $\Box$

\noindent{\bf Remark:} To obtain the result that the intrinsic tractor
parallel transport of $\nabla^{\ct_\Sigma}$ coincides with the
restriction of the ambient parallel transport of $\nabla^\cT$ to
sections of $\cT_\Sigma$ uses that $\Sigma$ is totally umbillic and
that $n^aC_{abcd}=0$ along $\Sigma$. These conditions are 
sufficient for the agreement of the connections. \quad \endrk

\subsection{Extending off $\Sigma$} \label{extS}
Given a conformal manifold $(\Sigma,[g_\Sigma])$ we may ask if this
can arise as the scale singularity set of a scalar negative almost
Einstein manifold. Narrowing the problem, we may begin with a fixed
smooth (or with specified regularity) codimension 1 embedding of $\S$
in a manifold $M$ and consider the Dirichlet-type problem of finding a
directed AE structure $(M,[g],I)$ with $(\Sigma,[g_\Sigma])$ as the
scale singularity set; issues include whether or not there is any
solution and, if there is, then whether $(\Sigma,[g_\Sigma]) $
determines $(M,[g],I)$ uniquely.  This is exactly the problem of
finding on $M$ a conformal structure $[g]$ and on this a solution
$\si$ to the conformally invariant equation $\nd_a\nd_b \si+P_{ab}
\si+\rho \bg_{ab}=0$ (i.e.\ \nn{primc}) such that $\S$ is the zero set
of $\si$ (and then there is the question of whether the pair
$([g],\si)$ is unique). We want to derive consequences of this
equation that make the nature of this problem more transparent. We
have seen already that this may be viewed as finding on $M$ a
conformal structure admitting a parallel tractor parallel tractor $I$
with $I|_\S$ agreeing with the normal tractor $N$ along $\S$.

The data on $\S$ is a conformal structure, and, for any solution
$[g_\S]$, is simply the pull back of the ambient conformal structure
$[g]$ on $M$.  By Theorem \ref{umbiltr} we know (at least to some
order along $\S$) how the ambient conformal curvature is related to
the intrinsic conformal curvature of $(\S,[g_\S])$. Thus it seems
natural to derive the equations which control how this extends off
$\S$. With less ambition we shall not attempt here to study the full
boundary problem. Rather we seek to find equations which control the
conformal curvature quantities off $\S$ and which are also well
defined along $\S$.

First note that it follows from the Bianchi identity \nn{bi2} that
Einstein manifolds $(M^d,\g)$ are Cotton, i.e. $A^{\g}_{abc}=0$. In
dimension $d=3$ the Weyl tensor vanishes identically and the Cotton
tensor is conformally invariant. Thus almost Einstein 3 manifolds are
Cotton and hence conformally flat. So if $(M^3,[g],I)$ is scalar
positive then it is a positive sectional curvature space form. If
$(M^3,[g],I)$ has $S(I)\leq 0$ then $I$ may have a scale singularity
set $\S$, but off this the structure $(M,\g)$ is either hyperbolic (if
$S(I)<0$) or locally Euclidean (if $S(I)=0$).

From \nn{divtrc} one easily concludes that on an Einstein manifold
 $(M^d,\g)$ the tractor curvature satisfies the (full) Yang-Mills
 equations, that is $\nabla^a \Omega_{ab}{}^C{}_D=0$ (see also
 \cite{GSS}) (where the connection $\nabla$ is in the scale $\g$). In
 dimension $d=4$ this equation is conformally invariant. Thus almost
 Einstein 4 manifolds are globally Yang-Mills. Combining with relevant results
 from Proposition \ref{aEsum} and Theorem \ref{umbiltr}, let us
 summarise .
\begin{proposition}\label{d4}
Let $(M^4,[g],I)$ be an almost Einstein manifold. Then the tractor curvature 
satisfies the conformally invariant Yang-Mills equations, 
$$
\nabla^a \Omega_{ab}{}^C{}_D=0. 
$$ If $I$ is scalar negative then along any singularity hypersurface
$\S$ of $I$ we have 
$$
C_{abcd}=0 \quad \mbox{and} \quad \Omega(u,v)=\Omega^\Sigma(u,v) \quad \mbox{ along } \Sigma
$$
where $u,v\in \Gamma(T\Sigma)$.
\end{proposition}
Note that in dimension 4 the tractor curvature is Yang Mills if and
only if the conformal structure is Bach flat. However the Proposition
suggests that it is useful to view the Bach flat condition as a
Yang-Mills equation in order to formulate an extension problem (or
boundary problem in the PE case). The additional data required
includes $n^a \Omega_{ab}{}^D{}_F$ along $\S$ which is equivalent to
$n^{b}A_{abc}|_\S$.  We note that in \cite{LeBrunH} LeBrun established
the existence and uniqueness of a real analytic self-dual
Poincar\'e-Einstein metric in dimension 4 defined near the
boundary with prescribed real analytic conformal infinity. If a
4-dimensional metric is self-dual then so is its tractor curvature and
hence the tractor connection is Yang-Mills.

Before we continue we need some further
notation.  Let us write $\# $ ({\em hash}) for the natural tensorial
action of sections $A$ of $\End(\cT)$ on tractor sections. For example,
on an covariant 2-tractor $T_{AB}$, we have
$$
A\hash T_{AB}=-A^{C}{}_AT_{CB}-A^{C}{}_BT_{AC}.
$$
If $A$ is skew for $\h$, then at each point, $A$ is $\frak{so}(
h)$-valued.  The hash action then commutes with the raising and
lowering of indices and preserves the $ SO(h)$-decomposition of
tractor bundles.  

As a section of the tensor square of the $\h$-skew bundle
endomorphisms of $\cT$, the curvature quantity $W$ has a double hash action
on tractors $T$; we write $W\hash \hash T$ for this.  Now for dimensions
$d\neq 4$ we use this to construct a Laplacian operator on
(possibly conformally weighted) tractor sections. For $T$ a section of 
$(\otimes^k \cT)[w]$  and $d\neq 4$ we make the definition
$$
\fb T : = (\Delta-wJ) T - \frac{1}{2(d-4)}W\hash\hash T .
$$
Then from this we obtain a variant of the usual tractor-D operator as follows:
$$
\fD_A T:=(d+2w-2)w Y_A T+ (d+2w-2)Z_{Aa}\nabla^a T 
+X_A \fb T. 
$$

In terms of this operator we have, 
\begin{theorem}\label{ext} Let $(M^d,[g],I)$ be an almost Einstein manifold.
Then if $d=4$ we have $W_{BCDE}=0$. In dimension 6 the
conformally invariant equation
$$
\fb W_{A_1A_2B_1B_2} = 0
$$
holds. 
In dimensions $d\neq 4$ we have 
\begin{equation}\label{extf}
I^A\fD_A W_{BCEF}=0.
\end{equation}
Also
\begin{equation}\label{WI0}
W_{BCEF}I^F=0= I^B W_{BCEF}.
\end{equation}
In particular if $I$ is scalar negative and $\S$ the singularity
hypersurface for $I$ then $W_{BCEF}N^F=0= N^B W_{BCEF}$ along $\S$.
\end{theorem}
\noindent{\bf Proof:} From \nn{Wform} it follows that, in dimension 4,
$W=0$ is equivalent setting the Bach tensor to zero and, as noted
earlier, this is equivalent to the conformal tractor connection being
Yang Mills. We have this, in particular, on almost Einstein manifolds.
Since $I$ is parallel it annihilates the tractor curvature, i.e.
$\Omega_{bc}{}^E{}_F I^E=0$. But since it is parallel and has conformal
weight 0, $I$ commutes with the tractor-D operator $\bD$.  It follows
from \nn{Wdef} that $W_{BCEF}I^F=0$. But $W$ has Weyl tensor type
symmetries, so \nn{WI0} follows. We also note here that since $I$
commutes with $\bD$, and, on the other hand, any contraction of $I$
with $W$ is zero, it follows by an elementary argument that $I$
commutes with $\fD$.

In dimension 6 we have from
\cite{GoPetobstrn} that
$$
\fb W_{A_1A_2B_1B_2} =
K X_{A_1}Z_{A_2}{}^aX_{B_1}Z_{B_2}{}^b \cB_{ab},
$$ where $K$ is a nonzero constant and $\cB_{ab}$ is the
Fefferman-Graham (obstruction) tensor (see also \cite{GoPrague}). The
sequentially labelled indices here are implicitly skewed over. But the
conformal invariant $\cB_{ab}$ is zero on Einstein manifolds
\cite{FGast,GrH,GoPetobstrn} and hence also (by continuity) on almost
Einstein manifolds.

It remains to establish \nn{extf}. Since $W$ has conformal weight
$-2$, it follows that when $d=6$ we have $I^A\fD_A W_{BCEF}= \si \fb
W_{A_1A_2B_1B_2}$, where as usual $\si$ denotes the conformal density
$X^AI_A$. Thus \nn{extf} holds in dimension 6.  Let us suppose now
that $d\neq 4,6$. Here we will use the link between the standard
tractor bundle on $(M,[g])$ and the Fefferman-Graham (FG) metric of
\cite{FGast,FGnew}.  This link was developed in
\cite{CapGoamb,GoPetCMP,BrGodeRham} but here we use especially the
notation and results from \cite{GoPetobstrn}. (It should be noted
however that here we use the opposite sign for the Laplacian.)  The
arguments we use below are a minor variation of similar developments
from those sources.

For a Riemannian conformal manifold $(M^d,[g])$ the ambient manifold
\cite{FGast} is a signature $(d+1,1)$ pseudo-Riemannian manifold with
$\cq$ as an embedded submanifold. There is some further background on
the FG metric in Section \ref{FGsec}. Suitably homogeneous tensor
fields on the ambient manifold, upon restriction to $\cq$, determine
tractor fields on the underlying conformal manifold.  In particular,
in dimensions other than 4, $W_{ABCD}$ is the tractor field equivalent
to $(d-4){\aR}_{ABCD}|_\cq$ where $\aR$ is the curvature of the FG
ambient metric. Under this correspondence the FG ambient metric
applied to tractors along $\cq$, descends to the tractor metric.
Ambient differential operators that are suitably tangential and
homogeneous (see e.g.\ \cite{BrGodeRham,GoPetobstrn}) also descend to
operators between tractor bundles or subquotients thereof. For example
the tractor connection arises from ambient parallel transport along
$\cq$.

On the FG ambient manifold let us define a Laplacian operator $\afl$ by
the formula
\[
\afl := \al - \frac{1}{2}\aR\#\# .
\]
Then in all dimensions $d\neq 4,6$, $\afl \aR|_\cq=0$, \cite[Section
3.2]{GoPetobstrn}. On the other hand $\fD_A$ corresponds to the
ambient operator
$$
(d+2w-2)\nda+\aX \afl=:\afD: \act^\Phi(w)\to \act\otimes  \act^\Phi(w-1)
$$ where $\act^\Phi(w)$ indicates the space of sections, homogeneous
of weight $w$, of some ambient tensor bundle. (NB: An ambient tensor $T$ is   homogeneous
of weight $w$ if and only if $\nda_{ \sX} T=w T$.)
From the Bianchi identity on the FG ambient manifold,  and 
the fact that $\afl \aR|_\cq=0$, it follows that 
 on the ambient manifold 
we have 
$$
\afD_{[A}\aR_{BC]DE}=0,
$$
along $\cq$. This descends to 
$$
\fD_{[A} W_{BC] DE}=0.
$$
So we have 
$$
I^A \fD_{[A} W_{BC] DE}=0.
$$
But since $I^A W_{ABDE}=0$ and $I$ commutes with $\fD$  \nn{extf} follows.
\quad $\Box$\\

From the Theorem we may conclude some restrictions on the intrinsic
conformal structure. For example we have the following. 
\begin{corollary}\label{Bflat}
Let $(M^5,[g],I)$ be a scalar negative almost Einstein manifold with
scale singularity set $\S\neq \emptyset$. Then the induced conformal structure 
$(\S^4,[g_\S])$ is Bach flat. 
\end{corollary}
\noindent{\bf Proof:} In dimensions $d\neq 6$ the equation \nn{extf},
i.e.\ $I^A\fD_A W=0$, on $M$ implies that along any scale singularity
subspace we $\S$ have 
$$\delta W =0
$$ where $\delta$ is the (conformally invariant) tractor twisted
conformal Robin operator \cite{BrGonon,GoMGC} applied to $W$; in terms
of $g$ we have $\delta W= n^a\nabla^g_a +2H^g W$ where $H^g$ is the
mean curvature of $\S$ and $\nabla$ is the usual (density coupled) the
tractor connection.

To simplify the presentation let us temporarily display the first two
abstract indices of the tractor $W$, but suppress the last pair; we
shall write $W_{BC}$ rather than $W_{BCDE}$. From the defining formula
\nn{Wdef} for $W_{BC}$ it follows easily that 
\begin{equation}\label{Ws}
W_{BC}= (d-4)Z_B{}^bZ_{C}{}^c \Omega_{bc}-X_BZ_C{}^b\nabla^a\Omega_{ab} +
X_CZ_B{}^b\nabla^a\Omega_{ab}.
\end{equation}
This is expression (13) from \cite{GoPetCMP}.  Exploiting Corollary
\ref{minc}, let us calculate in a metric $g$ with respect to which
$\S$ is totally geodesic.  Setting $d=5$,  applying
$\delta=n^a\nabla_a$ to \nn{Ws}, and using the tractor connection formulae
\nn{connids} we see that the coefficient of $Z_B{}^bZ_{C}{}^c$ is
$$
n^a\nabla_a\Omega_{bc}-n_b\nabla^a\Omega_{ac}+n_c\nabla^a\Omega_{ab} ,
$$ where $\Omega_{bc}$ is the tractor curvature of the ambient
conformal structure $(M,[g])$ (where we have suppressed the tractor
indices).  Evidently a part of the condition $\delta W|_\S=0$ is that
the last display is  zero along $\S$. Thus, in particular, 
$n^b$ contracted into this must vanish, that is
$$
n^bn^a\nabla_a\Omega_{bc}-\nabla^a\Omega_{ac}
+n_cn^b\nabla^a\Omega_{ab}=0\quad \mbox{along} \quad  \S.
$$ But using that $\S$ is totally geodesic and, from Theorem
\ref{umbiltr}, that $\Omega(u,v)=\Omega^\Sigma(u,v)$ along $\Sigma$
where $u,v\in \Gamma(T\Sigma)$, this exactly states that
$$
g^{ij}_\S \nabla^{\S}_i \Omega^\S_{jk}=0,
$$ where $g_\S$ is the intrinsic metric on $\S$ induced by $g$ and
$\Omega^\S_{jk}$ is the tractor curvature of its conformal class. Thus
the conformally invariant intrinsic tractor curvature of the
$(\S,g_\S)$ satisfies the Yang-Mills equations. As mention earlier, in
dimension 4 these are conformally invariant and are equivalent to the
structure $(\S,[g_\S])$ being Bach-flat.  \quad $\Box$

There is an analogue of this result for higher odd $d$, see Theorem
 \ref{FGflat} below. It is likely that there is a proof of Theorem
 \ref{FGflat} using only equation \nn{extf}, but certainly approaching
 this directly (as in the proof for $d=5$ above) would rapidly become
 technical for increasing dimension. Section \ref{FGsec} gives a
 simple and conceptual treatment, using the Fefferman-Graham metric.

\noindent{\bf Remarks:} From the equation \nn{extf} it follows that
the conformal aspects of the asymptotics of Poincar\'e-Einstein metrics are
controlled by the operator $I^A\fD_A$. 

In dimension 6 the main equation \nn{extf} (or equivalently \mbox{$\fb
W_{A_1A_2B_1B_2} = 0$}) is equivalent to requiring $(M,[g])$ to have
vanishing Fefferman-Graham tensor.  

In dimensions other than 3,4, and
6, and {\em off} $\S$, a key part of \nn{extf} is the harmonic
equation $\Delta C - \frac{1}{2} R\hash \hash C=0$ on the Weyl
curvature which holds on Cotton (and hence Einstein) manifolds, as
follows easily from the Bianchi identities \nn{WBianchiid} and
\nn{bi1}.  However in dimensions other than 3 we cannot conclude that
there is a scale for which an AE manifold is Cotton (everywhere). On
the other hand the equation \nn{extf} holds globally on an AE manifold
($d\neq 4$).  \quad \endrk

The following sheds some light on the meaning of the equation
\nn{extf} and its relation to possible boundary problems. This follows
easily from the Theorem and the definitions of the operators involved,
except we have also called on Corollary \ref{WW} below.
\begin{corollary}\label{interp}
On an Einstein manifold $(M^{(n+1)\geq 3}, \g)$ we have 
$$
(\Delta^{\g}+\frac{4J^{\g}}{n+1}(n-2))W - \frac{1}{2(n-3)}W\hash\hash W=0.
$$ In particular this holds on an almost manifold $(M^{(n+1)\geq
3},[g],I)$ off the zero set $\S$ of $\si=h(X,I)$.  If 
$|I|^2=1$  and $\S$ is non-empty, then on $M\setminus \S$ we have 
$$
(\Delta^{\g}-2(n-2))W - \frac{1}{2(n-3)}W\hash\hash W=0,
$$
while along the hypersurface $\S$ we have
$$
N\hook W=0, \quad \mbox{and, if~} n\geq 5, \quad (n-4)W|_\Sigma=(n-3)W^\Sigma,
$$
while,  if $n\neq 5$, 
$$
\delta W=0 \quad \mbox{along}\quad \S ,
$$ where $\delta $ is the conformal Robin operator applied to $W$; in
terms of $g$ we have $\delta W= n^a\nabla^g_a +2H^g W$ is the where
$H^g$ is the mean curvature of $\S$.
\end{corollary}
 It is shown in
\cite{GoMGC} that on densities $I^A D_A$ agrees with the Laplacians
arising in the scattering problems treated in \cite{GrZ}. The operator
$(\Delta^{\g}-2(n-2))$ here is a tractor twisted version of such. 
 We have used the $n=d-1$ here to simplify comparisons with
\cite{GoMGC} and \cite{GrZ}.

We have seen in dimension 3,4 and 6 that there are conformally
invariant equations controlling the conformal curvature of an AE
manifold. This is achieved trivially in dimension 3. As a final note
for this section we point out that there is an analogue of the results
for dimensions 4 and 6 to higher even dimensions. 
\begin{proposition}\label{conv}
Almost Einstein manifolds $(M^{d~{\rm even}},[g],I)$, $d\geq 4$,
satisfy the conformally invariant equation that the Fefferman-Graham
tensor vanishes.  This may be expressed in the form
\begin{equation}\label{FGe}
0=\fb_{d/2-2}W = \Delta^{d/2-2} W + \textit{lower~order~terms},
\end{equation}
where by $\fb_0$ and $\Delta^0$ we mean the operator given by
multiplication by 1.
\end{proposition}
The linear operator $\fb_{d/2-2}$ is constructed in
\cite{GoPetobstrn}, and the result here is an easy consequence of the
results there for Einstein manifolds. Once again on a scalar negative
AE manifold with a scale singularity set $\S$, \nn{FGe} expresses the
vanishing Fefferman-Graham tensor condition in a form suitable to link
with the conformal curvature data on $\S$ (using Corollary \ref{WW},
or for $n=3,4$ Theorem \ref{umbiltr}). It should be interesting to
construct compatible conformal boundary operators for $W$ along
embedded submanifolds $\S$ so that these yield a well posed and
conformal elliptic problem for the conformal curvature $\Omega$. 
Close analogues of the conformal boundary operators developed in
\cite{BrGonon} should play a role.\\
\noindent{\bf Remark:} Note that conformal equations, such as
\nn{FGe}, offer the chance to split the problem of finding almost
Einstein structures (or Poincar\'e-Einstein metrics) into a conformal
problem, say controlled by \nn{FGe} with further boundary operators
along $\S$, and a second part where one would find a compatible
``scale'' $\si$. We should expect that a solution to the conformal
problem is necessary but in general not sufficient.  However one may
ask if (in Riemannian signature and say on closed even manifolds) \nn{FGe} plus
the (clearly necessary) vanishing of the conformal invariant
$$ 
\Omega_{ab}{}^C{}_{F_1} \Omega_{cd}{}^D{}_{F_2} \cdots
\Omega_{ef}{}^E{}_{F_{d+1}},
$$ where the sequentially labelled indices are skewed over, is
sufficient for a conformal manifold to necessarily admit an almost
Einstein structure locally. A corresponding global question is whether
a smooth section $K$ of $\cT$ satisfying $\Omega_{ab}{}^C{}_DK^D=0$
plus \nn{FGe} is sufficient to conclude that the  conformal structure
on a closed even manifold admits a directed almost Einstein
structure. In dimension 4 there is a positive answer to this if we
restrict to $K$ such that $h(X,K)$ is non-vanishing \cite{GoNagy}; in this
case the structure must be Einstein.\quad \endrk

\section{Examples and the model} \label{exsS}

\subsection{The model -- almost Einstein structures on the sphere}
\label{model}

\begin{proposition}\label{modp}
  The $d$-sphere, with its standard conformal structure, admits a
  $(d+2)$-dimensional space of compatible directed almost Einstein
  structures.  For each $S\in \mathbb{R}$ there is an almost Einstein
  structure $I$ on $\bS^d$ with $S(I)=S$.
\end{proposition}

\noindent The AE structures on the sphere also may be viewed as
examples of ASC structures on the sphere. In any case we shall see
that, in a sense, ``most'' of these are scalar negative (which might
at first seem counterintuitive).

Before we prove this let us recall the construction of the standard
conformal structure on the sphere.  Consider a $(d+2)$-dimensional
real vector space $\bV$ equipped with a non-degenerate bilinear form
$\mathcal{H}$ of signature $(d+1,1)$. The {\em null cone} $\cN$ of
zero-length vectors form a quadratic variety.  Let us write $\cN_+$
for the forward part of $\cN\setminus \{0 \}$.  Under the ray
projectivisation of $\bV$ the forward cone $\cN_+$ is mapped to a
quadric in $\P_+(\bV) \cong \bS^{d+1}$.  This image is topologically a
sphere $\bS^d$ and we will write $\pi$ for the submersion $\cN_+\to
\bS^d$. Each point $p\in \cN_+$ determines a
positive definite inner product on $T_{x=\pi{p}}\bS^d$ by
$g_x(u,v)=\mathcal{H}_p(u',v')$ where $u',v'\in T_p\cN_+$ are lifts of
$u,v\in T_x\bS^d$. For a given vector $u\in T_x \bS^{d}$ two lifts to
$p\in \cN_+$ differ by a vertical vector field. Since any vertical
vector is normal (with respect to $\mathcal{H}$) to the cone it
follows that $g_x$ is independent of the choices of lifts. Clearly
then, each section of $\pi$ determines a metric on $\bS$ and by
construction this is smooth if the section is. (Evidently the metric
agrees with the pull-back of $\mathcal{H}$ via the section concerned.) Now,
viewed as a metric on $T\mathbb{R}^{d+2}$, $\mathcal{H}$ is
homogeneous of degree 2 with respect to the standard Euler vector
field $E$ on $\bV$, that is $\cL_E \mathcal{H}=2 \mathcal{H}$,
where $\cL$ denotes the Lie derivative. In particular this holds on
the cone, which we note is generated by $E$.

Write $\bg$ for the restriction of $\mathcal{H}$ to vector fields in
$T\cN_+$ which are the lifts of vector fields on $\bS^d$. Then for any
pair $u,v\in \Gamma(T\bS^d)$, with lifts to vector fields $u',v'$ on
$\cN_+$, $\bg(u',v')$ is a function on $\cN_+$ homogeneous of degree
2, and which is independent of how the vector fields were
lifted. Evidently $\cN_+$ may be identified with the total space of a
bundle of conformally related metrics.  Thus $\bg(u',v')$ may be
identified with a conformal density of weight $2$ on $\bS^d$. That is,
this construction determines a section of $S^2T^*\bS^d\otimes E[2]$
that we shall also denote by $\bg$. By construction this is a
conformal metric (see Section \ref{al}) on $\bS^d$.  Fix a future
pointing vector $I$ in $\bV$ with $|I|^2:=\mathcal{H}(I,I)=-1$.
Regarding $\bV$ as an affine space, view $I$ as a constant section of
$T\bV$. Write $X^A$ for standard coordinates on $\bV$ (i.e. via an
isomorphism $\bV\cong \mathbb{R}^{d+2}$). It is straightforward to
verify that the $\mathcal{H}(I,X)=1$ hyperplane meets $\cN_+$ in a
copy of $\bS^d$ and the metric induced by this section of $\pi$ is the
standard metric on $\bS^d$. Thus $\bg$ is a standard conformal
structure on the sphere. We are ready to prove the Proposition.

\smallskip 

\noindent{\bf Proof of Proposition \ref{modp}:} It is easily verified
that $G:=SO(\mathcal{H})\cong SO_0(d+1,1)$ (the identity connected
component of the Lorentz group) acts transitively on the sphere. Thus
the conformal sphere may be identified with $G/P$ where $P$ is the
parabolic subgroup of $G$ which stabilises a nominated ray in $\cN_+$.
Now $G\to G/P$ may be viewed as a flat Cartan bundle over $G/P=\bS^d$
and the standard tractor bundle $\cT$ is $G\times _P \mathbb{V}$ where
$\mathbb{V}$ is viewed as a $P$-module, by restriction. Here $G\times
_P \mathbb{V}=G\times \bV/\sim $ where the equivalence relation is
$(gp,v)\sim (g,p\cdot v)$ with $g\in G$, $p\in P$ and where ``$\cdot
$'' indicates the standard representation of $G$ on $\bV$.  The bundle
$G\times _P \mathbb{V}$ is trivialised canonically by the map
$(g,v)\mapsto (gP,g\cdot v)$ and so we have a connection $\nabla^\cT$
on $\cT$ induced from the trivial connection on $(G/P)\times \bV$. It
is straightforward to verify that this is the normal tractor
connection. (In fact this is essentially a tautology; one view the
idea of a normal conformal connection tractor as modelled on this
homogeneous case.)  Thus in this case the tractor connection is
globally flat, with the bundle $\cT$ admits $(d+2)$ linearly
independent parallel sections.  \quad $\Box$

Using the embedding of $\cN_+$ in $\bV$ we can explicitly describe the
almost Einstein structures of the Proposition.  For example we may
construct a scalar negative AE structure on $\bS^d$ as follows.  Take
a vector $I\in \bV$ of length 1 (i.e. $|I|^2=1$).  We shall use the
same notation for the covector $\mathcal{H}(I,\cdot)$. By the standard
parallel transport (of $\bV$ viewed as an affine structure) view this
as a constant section of $T^*\bV$.  Then as above, writing $X^A$ for
standard coordinates on $\bV$, the intersection of the hyperplane
$I_AX^A=1$ with $\cN_+$, which we shall denote $S_+$, is a section of
$\pi$ over an open cap $C_+$ of the sphere.  Similarly the
intersection of the hyperplane $I_AX^A=-1$ with $\cN_+$, which we
shall denote $S_-$, is a section of $\pi$ over another open cap $C_-$
of the sphere.  On the other hand the hyperplane $I_AX^A=0$ (parallel
to the previous) intersects $\cN_+$ in a cone of one lower
dimension. The image $\S$ of this under $\pi$ is a copy of $\bS^{n}$
embedded in $\bS^d$ (where as usual $d=n+1$). It is easily deduced
that $\bS^d$ is the union of the three submanifolds $C_-$, $\S$, and
$C_+$ and that, for example, with respect to (a restriction of) the
smooth structure structure on $\bS^d$, the embedded $\S$ is a boundary
for its union with $C_+$.
This follows because any forward null ray though the origin and
parallel to the $I_AX^A=1$ hyperplane lies in the hyperplane
$I_AX^A=0$, whereas every other forward null ray through the original
meets either the $I_AX^A=1$ hyperplane or the $I_AX^A=-1$
hyperplane. 
Let us write $\g$ for the metric that the sections $S_\pm$  give on
$C_\pm$. 
Note that the hypersurface $\bS^{n}$ canonically has
no more than a conformal structure. This may obviously be viewed as
arising as a restriction of the conformal structure on
$\bS^{d}$. Equivalently we may view its conformal structure as arising
in the same way as the conformal structure on $\bS^d$, except in this
case by the restriction of $\pi$ to the sub-cone $I_AX^A=0$ in
$\cN_+$, and from (the restriction of) $\mathcal{H}$ along this sub-cone. 
In the following we write $g$ to denote any
metric from the standard conformal class on $\bS^d$.
Note that on $C_\pm$ this is
conformally related to $\g$.

Now let us henceforth identify, without further mention, each function
on $\cN_+$ which is homogeneous of degree $w\in \R$ with the
corresponding conformal density of weight $w$. With $\si:=I_AX^A$, as
above, note that $\si^{-2}\bg$ is homogeneous of degree 0 on $\cN_+$
and agrees with the restriction of $\mathcal{H}$ along $S_\pm$. Thus
on $C_\pm$ we have $\si^{-2}\bg=\g$; $\si^{-2}\bg$ recovers the metric
determined by $S_\pm$.  Similarly on $\bS^d$ we have $g=\tau^{-2}\bg$,
where $\tau$ is a non-vanishing conformal density of weight 1.  So on
$C_+\cup C_-$, $\g=s^{-2}g$ where $s$ is the {\em function}
$\si/\tau$. We see that $\g$ is conformally compact on $\bS^d\setminus
C_+$, and also on $\bS^d\setminus C_-$.

We may now understand this structure via the tractor bundle on
$\bS^d$. Let us write $\rho^t$ for the natural action of $\R_+$ on
$\cN_+$ and then $\rho^t_*$ for the derivative of this. Now modify the
latter action on $T\bV$ by rescaling: we write $t^{-1}\rho^t_*$ for
the action of $\R_+$ on $T\bV$ which takes $u\in T_p\bV$ to
$t^{-1}(\rho_*^t u)\in T_{\rho^t(p)}\bV$. Note that $u$ and
$t^{-1}(\rho_*^t u) $ are mutually parallel, according to the affine structure
on $\bV$. It is easily verified that the quotient of $T\bV|_{\cN_+}$ by
the $ \R_+$ action just defined is a rank $d+2$ vector bundle $\cT$ on
$M$. Obviously the parallel transport of $\bV$ determines a parallel
transport on $\cT$, that is a connection $\nabla^\cT$. Since $\bV$ is
totally parallel this connection is flat. The twisting of $\rho^t_*$
to $t^{-1}\rho^t_*$ is designed so that the metric $\mathcal{H}$ on
$\R^{d+2}$ also descends to give a (signature $(d+1,1)$) metric $h$ on
$\cT$ and clearly this is preserved by the connection.  In fact
$(\ct,h,\nabla^\cT)$ is the usual normal standard tractor bundle. This
is proved under far more general circumstances in \cite{CapGoamb} (see
also \cite{GoPetCMP}); it is shown there that the tractor bundle may
be recovered from the Fefferman-Graham ambient metric by an argument
generalising that above. In this picture the Euler vector field
$E=X^A\partial/\partial X^A$ (using the summation convention), which
generates the fibres of $\pi$, descends to the canonical tractor field
$X\in \ct[1]$.

It follows from these observations that, since the vector field $I$ is
parallel on $\bV$, its restriction to $\cN_+$ is equivalent to a
parallel section of $\cT$; we shall also denote this by $I$. So this
is an almost Einstein structure on $\bS^d$; $|I|^2=1$ means that the
almost Einstein structure we recover has $S(\si)=-1$, whence has
$\Ric(\g)=-n \g$ on $C_\pm$.  The zero set for $\si=h(X,I)$ is exactly
$\S$. So we see that $(\bS^d,[g],I)$ is an almost Einstein
manifold. Since it is conformally flat with $S(\si)=-1$ it is what may
be termed an {\em almost hyperbolic} structure on the sphere.  The
fact that along $\S$ the parallel tractor $I$ gives the normal tractor
$N$ is especially natural in this picture since $\S$ is determined by
a hyperplane orthogonal to $I$. Finally we observe that it follows
from Proposition \ref{peae} that the spaces $(\bS^d\setminus C_\pm,[g],I)$
are Poincar\'e-Einstein manifolds, in fact each equivalent to the
conformal compactification of the hyperbolic ball.

Since the group $G$ acts transitively on length 1 spacelike vectors,
from the picture above we see that any scalar negative AE structure on
the sphere is related to the one constructed by a conformal
transformation after an $\mathbb{R_+}$ action on the parallel tractor
$I$. 

The scalar flat almost Einstein structures are obtained by a similar
construction to the scalar negative case above. Note that if $I$ is a
non-zero null vector in $\bV$ then the hyperplane $\mathcal{H}(I,X)=1$
meets all future null rays in $\cN_+$ except the one parallel to $I$.
So the almost Einstein structure determined by $I$ has a single
isolated point of scale singularity.  The Einstein metric $\g$ is
conformally related to the round metric, and $|I|^2=0$ means that
$S(I)=0$ and so $\g$ is flat; this is the usual Euclidean structure on
the sphere minus a point. It is straightforward to conclude that the
map, along null generators, relating this Euclidean almost Einstein
structure and the standard sphere embedded in the cone (as described
earlier) is the usual stereographic projection.

In a partial summary then, if $I_1$ and $I_2$ are constant vectors in
$\bV$ with $|I_1|^2=-1$ and $|I_2|^2=1$ then, as parallel tractors on
$\bS^d$ these determine, respectively the standard sphere metric and
almost hyperbolic structures.  We can interpolate between these via
Corollary \ref{betw} and we note that for some $t\in \mathbb{R}$ the
parallel tractor $I_t:=(\sin t)I_1+ (\cos t)I_2$ is null and so
determines a Ricci flat structure in the conformal class, that is a
Euclidean metric on the sphere minus a point. For each $t\in
\mathbb{R}$ the isotropy subgroup $G_{I_t}$ of $G=SO_0(\mathcal{H})$
fixing the vector $I_t$ clearly acts transitively and by isometries on
the connected components of $\bS^d\setminus \S_t$, where $\S_t$ is the
scale singularity set of $I_t$.

\subsection{Doubling and almost hyperbolic constructions} \label{dou}
One route to 
constructing further compact almost Einstein manifolds is via the
doubling of compact Poincar\'e-Einstein manifolds. So suppose that $M$
is a compact Poincar\'e-Einstein with conformal infinity $\Sigma$. The
double we seek is a gluing along $\Sigma$,
$$
M_{(2)}:= (M\sqcup M)/\Sigma
$$
where the identification of the two copies of $\Sigma$ is the obvious
one. As pointed out in \cite{MM}, for example, this may be equipped
with a smooth structure compatible with the smooth structure on $M$
and so that the natural involution exchanging the factors is also
smooth. Now extend the PE metric $\g$ of $M$ to a metric on $M_{(2)}$
by symmetry. This will be smooth if $\g$ is {\em even} in the sense
of \cite[Section 4]{FGnew} (following \cite{GL}):  Locally along the
collar the metric may be put in normal form, relative to some
$g^\Sigma$ from the conformal class on $\Sigma$,
\begin{equation}\label{normalf}
\g =s^{-2}(ds^2+ g^\Sigma_s) 
\end{equation}
where $s$ satisfies $|ds|^2_{g}=1$ and $g_s$ is a 1-parameter family
of metrics on $\Sigma$ such that $g^\Sigma_0=g^\Sigma$. The metric is
{\em even} if for each point of $\Sigma$, and with the metric $\g$ in
this form, we have that $ds^2+ g_s$ is the restriction to $M\times
[0,\infty)$ of a smooth metric $g$ on a neighbourhood
$\mathcal{U}\subset \Sigma\times (-\infty,\infty)$ such that $\mathcal{U}$ and
$g$ are invariant under the map $s\mapsto -s$.

Infinite volume hyperbolic manifolds provide a source of
even PE manifolds. From Theorem 7.4 in \cite{FGnew} (building on
\cite{E1,SS}) we have that if $(M,\g)$ is a hyperbolic PE manifold
then locally along the conformal infinity $\Sigma$ it may be put in
the normal form \nn{normalf} where (in terms of local coordinates
$(s,x^i)$, with $x^i$ the coordinates on $\Sigma$) we have
$$
(g_s)_{ij}=g^\Sigma_{ij}-P^\Sigma_{ij}s^2 
+\frac{1}{4}g^{kl}_\Sigma P^\Sigma_{ik}P^\Sigma_{lj} s^4;
$$ here if $d\geq 4$ then $P^\Sigma_{ij}$ is the intrinsic Schouten
tensor of $g^\Sigma$, while if $d=3$ then $P^\Sigma_{ij}$ is a
symmetric 2-tensor on $\Sigma$ satisfying
$2g_\Sigma^{ij}P^\Sigma_{ij}=\Sc^{g^\Sigma}$ and $2
g_\Sigma^{ij}\nabla_i P^\Sigma_{jk} = \nabla^\Sigma_k
\Sc^{g^\Sigma}$. In this case $\g$ is manifestly even.

Let $\Gamma$ be a convex co-compact,
torsion-free, discrete group of orientation preserving isometries of
$\mathbb{H}^d$. Then the orbit space $M_+:=\Gamma \backslash \mathbb{H}^d$
  is a hyperbolic manifold of
infinite volume. Such $M_+$ may be conformally compactified
\cite{M-hodge,MM} to yield a (hyperbolic) PE manifold. Thus, by the
doubling construction, to each group $\Gamma$ as above we may
associate a closed almost Einstein structure.
 
Rather than the usual model of the hyperbolic ball we may realise
$\mathbb{H}^d$ as a hyperbolic cap of the sphere as described in
Section \ref{model}; it is not difficult to see that we may arrange
that the cap is the right-hemisphere of a standard round sphere, where
the latter is given also as a section of the cone as Section
\ref{model}.  Then the convex co-compact $\Gamma$ arises as a discrete
subgroup of $G_I\subset G=SO(\mathcal{H})$, where $G_I$ is the
isotropy subgroup of $G$ which fixes the length 1 parallel tractor $I$
defining the hyperbolic manifold.  Now $\Gamma$ also acts on the
hyperbolic left-hemisphere and, by symmetry, in both cases the
conformal infinity $\Sigma_\Gamma$ may be identified with the orbit
space $\Gamma \backslash \Omega_\Gamma (\S)$ where $\Omega_\Gamma
(\S)$ is the open subset of the sphere $\S\cong \bS^n$ where $\Gamma$
acts properly discontinuously.  The smooth structure on the doubling
of $M$ is the usual smooth structure on the sphere and the action of
$\Gamma$ on $\bS^d$ evidently preserves the solution of \nn{primc}
giving $I$.  In this sense we may view the doubling of $M$ as arising
from the orbit space of $\Gamma$ on $\bS^d$ equipped with a standard
almost hyperbolic structure (i.e.\ hyperbolic almost Einstein
structure as in Section \ref{model}), but where we first remove the
limit of this action in $\S$. The AE structure on this is a solution
to \nn{primc} descended from a solution on $\bS^d$.

\section{The Fefferman-Graham metric for an AE manifold and obstructions}\label{FGsec}

The Fefferman-Graham tensor (also called ``the obstruction tensor'')
is a natural conformally invariant symmetric trace-free 2-tensor
$\mathcal{B}_{ab}$ on manifolds of even dimension $n$ that has the
form $\Delta^{n/2-2} \nd^c\nd^d C_{acbd} +lower~order~terms$. In the
case $n=4$ it agrees with the Bach tensor while in higher even
dimensions it is due to Fefferman and Graham \cite{FGast}. In
Corollary \ref{Bflat} we found that the Bach tensor necessarily
vanishes on the scale singularity set of AE 5-manifolds. Here we prove
the analogue of that result for higher odd dimensions. In the process
of proving this we obtain an extension to Theorem \ref{umbiltr}.

For $\pi:\cq\to M^d$ a Riemannian conformal structure, let us use
$\rho $ to denote the ${\Bbb R}_+$ action on $ \cq$ given by $\rho(s)
(x,g_x)=(x,s^2g_x)$.  An {\em ambient manifold\/} is a smooth
$(d+2)$-manifold $\aM$ endowed with a free $\Bbb R_+$--action $\rho$
and an $\Bbb R_+$--equivariant embedding $i:\cq\to\aM$.  We write
$\X\in\Gamma(T \aM)$ for the fundamental field generating the $\Bbb
R_+$--action.  That is, for $f\in C^\infty(\aM)$ and $ u\in \aM$, we
have $\X f(u)=(d/dt)f(\rho(e^t)u)|_{t=0}$.  For an ambient manifold
$\aM$, an {\em ambient metric\/} is a pseudo--Riemannian metric $\h$
of signature $(d+1,1)$ on $\aM$ satisfying the conditions: (i) $\Cal
L_{\sX}\h=2\h$, where $\Cal L_{\sX}$ denotes the Lie derivative by
$\X$; (ii) for $u=(x,g_x)\in \cq$ and $\xi,\eta\in T_u\cq$, we have
$\h(i_*\xi,i_*\eta)=g_x(\pi_*\xi,\pi_*\eta)$. In \cite{FGast} (and see
\cite{FGnew}) Fefferman and Graham considered formally the Gursat
problem of obtaining $\Ric(\h)=0$. They proved that for the case of
$d=2$ and $d\geq 3$ odd this may be achieved to all orders, while for
$d\geq 4$ even the problem is obstructed at finite order by the tensor
$\mathcal{B}_{ab}$; for $d$ even one may obtain $\Ric(\h)=0$ up to the
addition of terms vanishing to order $d/2-1$. (See \cite{FGnew} for
the statements concerning uniqueness. For extracting results via
tractors we do not need this, as discussed in e.g.\
\cite{CapGoamb,GoPetCMP}.) We shall henceforth call any (approximately
or otherwise) Ricci-flat ambient metric a {\em Fefferman-Graham
metric}.

Since an AE manifold $(M^d,[g],I)$ has, by definition, a conformal
structure we may construct the Fefferman-Graham metric, as for any
conformal manifold.  We have already exploited this in the proof of
Theorem \ref{ext}. On the other hand if $S(I)<0$ and the scale
singularity set $\S$ is non-empty then, as discussed in Section
\ref{MvS}, this embedded $n$-manifold ($n=d-1$) has induced on it a conformal
structure $(\S,[g_\S])$. We may ask how the Fefferman-Graham metric
for $(\S,[g_\S])$ is related to the Fefferman-Graham metric for
$(M^d,[g])$. In Theorem \ref{amb} below for $n$ even (so $d=n+1$ odd)
we show that $(\S,[g_\S])$ admits a Fefferman-Graham which is
Ricci-flat to all orders.  Thus we obtain
the following.
\begin{theorem} \label{FGflat} Suppose than $(\S^n,[g_\S])$ is the scale 
singularity space of a scalar negative almost Einstein manifold, then the
Fefferman-Graham tensor of $(\S^n,[g_\S])$ is zero.
\end{theorem}
\noindent Using Theorem \ref{maina}, this result also follows from
 \cite[Theorem 4.8]{FGnew} or \cite[Theorem 2.1]{GrH}. The proof here
 follows a rather different tack.

In the subsequent discussion of ambient metrics all results can be
assumed to hold formally to all orders unless stated otherwise.  We
typically use bold symbols or tilded symbols for the objects on
$\tilde{M}$. For example $\boldsymbol{\nabla}$ is the Levi-Civita
connection on $\tilde{M}$.
It is assumed the reader is somewhat familiar with
treatments of Fefferman-Graham metrics. In particular, as used in the
proof of Theorem \ref{ext}, we use that that suitably homogeneous
tensor fields of $\tilde{M}|_\cq$ correspond to tractor fields.  The
notation and approach here follows that in
\cite{BrGodeRham,CapGoamb,GoPetobstrn}.

\begin{lemma}\label{silift} Let $(M^d,[g],I)$ be an AE manifold with 
$d\geq 3$ odd.  There is a parallel 1-form field $\boldsymbol{I}$ on
$\tilde{M}$ such that $\boldsymbol{I}|_\cq$ is the homogeneous (of
weight 0) section of $T^*\tilde{M}|_\cq$ corresponding to $I$.
\end{lemma}
\noindent{\bf Proof:} Let $\si:=h(X,I)$, as usual. This corresponds to
function on $\cq$ homogeneous of degree 1. Since $d$ is odd, this may
be extended ``harmonically'' to all orders (e.g.\ \cite{GJMS}). That
is there is a homogeneous degree 1 function $\bs{\si}$ on $\tilde{M}$
such that $\bs{\Delta} \bs{\si}=0$ and $\bs{\si}|_\cq$ is the
homogeneous function corresponding to the conformal density $\si$.

The operator $\tilde{\bD}_A=(d+2w-2)\bs{\nd}_A+\bs{X}_A\bs{\Delta}$,
on $\tilde{M}$, corresponds to the tractor-$D$ operator $\bD$
\cite{CapGoamb,GoPetCMP}. This acts tangentially along $\cq$ in the
sense of \cite{BrGodeRham} and \cite{GoPetobstrn}. Define
$\bs{I}_A:=\bs{\nabla}_A\bs{\si}=\frac{1}{d}\tilde{\bD}_A
\bs{\si}$. This has the required properties.  Obviously
$\boldsymbol{I}|_\cq$ corresponds to $I_A=\frac{1}{d}\bD_A
\si$. (Recall on a density $\si$ of weight 1, $\frac{1}{d}\bD_A \si= D
\si$.)  Now note that
$\bs{\Delta}\bs{I}_B=\bs{\Delta}\bs{\nabla}_B\si=
\bs{\nabla}_B\bs{\Delta}\si = 0$, to all orders, as the FG metric is
Ricci flat to all orders. Now $\bD_A I_B=0$ on $M$, and so
$\tilde{\bD}_A \bs{I}_B|_\cq=0$. Using the previous result we conclude
$((d-2)\bs{\nabla}_A \bs{I}_B)|_\cq=0$.  Now by induction we get that
$\bs{I}_A$ is parallel to all orders: Suppose that
$\bs{\nabla}_{A_2}\cdots \bs{\nabla}_{A_{i+1}}\bs{I}_B|_\cq=0$ for
$i=1,\cdots ,k$ then, since $\tilde{\bD}$ acts tangentially, we get
$$
\tilde{\bD}_{A_1} \bs{\nabla}_{A_2}\cdots \bs{\nabla}_{A_{k+1}}\bs{I}_B|_\cq=0.
$$
Thus, along $\cq$,
$$ (d-2k-2)\bs{\nabla}_{A_1}\bs{\nabla}_{A_2}\cdots
\bs{\nabla}_{A_{k+1}}\bs{I}_B + \bs{X}_{A_1} \bs{\Delta}
\bs{\nabla}_{A_2}\cdots \bs{\nabla}_{A_{k+1}}\bs{I}_B =0
$$ To study the second term we may commute the Laplacian $\bs{\Delta}$
to the right of the $\bs{\nabla}$'s with free indices. We see then that
this entire term drops out as $\bs{\Delta}\bs{I}$ vanishes to all orders
while the other terms pick up curvature and hence involve at most
$(k-1)$ derivatives of $\bs{I}$ (and so exit by the inductive
hypothesis). On the other hand, since $d$ is odd, $(d-2k-2)\neq 0 $.
\quad $\Box$

Using this we obtain the key result.
\begin{theorem}\label{amb}
Let $(M^d,[g],I)$ be a scalar negative AE manifold with $d\geq 3$ odd and $\S\neq \emptyset$. Write 
 $\bs{I}$ for the parallel 1-form field on $\tilde{M}$ corresponding
(as in the Lemma above) to $I$.  Write $\bs\S$ for the hypersurface
given as the zero set of $\bs{\si}:=\bs{h}(\bs{X},\bs{I})$. This has a
metric $\bs{h}_{\bs{\S}}$ induced from $\bs{h}$, it is totally
geodesic, and $(\bs{\S}, \bs{h}_{\bs{\S}})$ is a Fefferman-Graham
metric for $(\S,[g_\S])$, which is formally smooth and Ricci-flat to
all orders.
\end{theorem}
\noindent{\bf Proof:}
First some observations. It is clear that
$\bs{\S}$ is a smooth hypersurface and its intersection with $\cq$ is
the inverse image of $\S$ with respect to the standard map $\cq\to
M$. Since the conformal structure $[g_\S]$ of $\S$ is induced from the
conformal structure of the ambient space $(M,[g])$ it follows easily
that restricted to the tangents of this intersection
$\bs{h}_{\bs{\S}}$ agrees with the tautological 2-form (which we have
since the intersection of $\bs{\S}$ with $\cq$ is naturally identified
with the bundle of metrics in the conformal class over $(\S,[g_\S])$.
 
 Since $\bs{\si}$ is
homogeneous of degree 1 we have $\mathcal{L}_{\bs{X}}\bs{\si}=\bs{\si}$ and so along
$\bs{\S}$ the field $\bs{X}$ is everywhere tangent to $\bs{\S}$. Clearly
$\mathcal{L}_{\bs{X}}\bs{h}^{\bs{\S}}= 2 \bs{h}^{\bs{\S}} $ from the
analogous property for $\bs{h}$.

 Since $\bs{d}\bs{\si}$ is parallel $h^{-1}(\bs{d}\bs{\si}, \bs{d}\bs{\si})$
is constant and agrees with $|I|^2>0$. In particular, along $\bs{\S}$,
$\bs{N}:=\bs{d}\bs{\si}$ gives a parallel conormal field for $\bs{\S}$, of
non-zero pointwise length. So $\bs{\S}$ is totally geodesic.

Once again using that $\bs{d}\bs{\si}$ is parallel we  have that
$\bs{R}_{AB}{}^C{}_D \bs{N}_C=0$ along $\bs{\S}$. It follows that the
intrinsic Ricci curvature $\bs{\Ric}^{\bs{\S}}$ agrees with the
 tangential restriction of the ambient Ricci curvature. But the
latter is everywhere zero to all orders, and therefore so is
$\bs{\Ric}^{\bs{\S}}$. \quad $\Box$

\begin{corollary} \label{WW}If $(M^d,[g],I)$ is a scalar negative AE structure 
with $d\geq 6$ and a non-empty scale singularity space $\S$, 
then $(d-5)W|_\S= (d-4)W^\S$.
\end{corollary}
\noindent In the Corollary we view, by trivial extension, $W^\S$ as a section of $\otimes^4 \cT$.
\noindent{\bf Proof:} If $d\geq 7$ is odd then this immediate from the
proof above. Since $\bs{\S}$ is totally geodesic and $\bs{N}$
annihilates the curvature $\bs{R}$ of the Fefferman-Graham metric
$\bs{h}$, it follows that along $\bs{\S}$ we have
$\bs{R}=\bs{R}^{\bs{\S}}$ (using a trivial extension to view
$\bs{R}^{\bs{\S}}$ as a section of $\otimes T^*\tilde{M}$).  But, as
used in Section \ref{extS}, $(d-4)\bs{R}|_\cq$ is the ambient tensor
field field equivalent to the tractor $W$; $(n-4)\bs{R}^{\bs{\S}}$ 
similarly corresponds to $W^\S$.

If $d\geq 6$ is even it is straightforward to verify that 
the results in Lemma \ref{silift} and in Theorem \ref{amb} hold to sufficient order to obtain the result here.
\quad $\Box$ \\ In some sense the Corollary applies to all dimensions
$d\geq 3$ except for $d=4$ as follows.  When $d=5$ since $[g_\S]$ is Bach
flat we have $W^\S=0$. However $W|_\S$ gives a tractor field
equivalent to $\bs{R}^{\S}$. For $d=3$ AE manifolds both $\bs{R}$ and
$\bs{R}^{\bs{\S}}$ are zero.

\end{document}